%% file: ex_article.tex
\pdfoutput=1
\documentclass[hidelinks,onefignum,onetabnum]{siamart220329}

\input{ex_shared}

\ifpdf
\hypersetup{
  pdftitle={High-performance matrix-free unfitted finite element operator evaluation},
  pdfauthor={M. Bergbauer, P. Munch, W. A. Wall, and M. Kronbichler}
}
\fi

\begin{document}

\maketitle

\begin{abstract}
  Unfitted finite element methods, like CutFEM, have traditionally been implemented in a matrix-based fashion, where a sparse matrix is assembled and later applied to vectors while solving the resulting linear system. With the goal of increasing performance and enabling algorithms with polynomial spaces of higher degrees, this contribution chooses a more abstract approach by matrix-free evaluation of the operator action on vectors instead. The proposed method loops over cells and locally evaluates the cell, face, and interface integrals, including the contributions from cut cells and the different means of stabilization.
  The main challenge is the efficient numerical evaluation of terms in the weak form with unstructured quadrature points arising from the unfitted discretization in cells cut by the interface. We present design choices and performance optimizations for tensor-product elements and demonstrate the performance by means of benchmarks and application examples. We demonstrate a speedup of more than one order of magnitude for the operator evaluation of a discontinuous Galerkin discretization with polynomial degree three compared to a sparse matrix-vector product and develop performance models to quantify the performance properties over a wide range of polynomial degrees.
\end{abstract}

\begin{keywords}
matrix-free, high-order, high-performance, tensor-product, unfitted finite element method
\end{keywords}

\begin{MSCcodes}
  65N30, 65Y20, 65Y05, 68W10
\end{MSCcodes}

\section{Introduction}

For solving the large sparse linear systems resulting from finite-element
discretizations of PDEs, matrix-free evaluation of
the operator action
has become increasingly popular against the traditional sparse matrix-based
realizations during the past decade. Besides increasing the maximal feasible problem size
via reduced memory requirements, the primary benefit is the
higher throughput of the operator evaluation when used within iterative
solvers~\cite{Kronbichler2018}. Matrix-based iterative solvers are usually bound by
memory bandwidth, a result of the low arithmetic intensity of sparse matrix-vector products and the
arithmetic capabilities of modern hardware, which have been growing more
rapidly than the available memory bandwidth, both on CPU and GPU
architectures, for several decades already. Matrix-free methods align better with the
underlying hardware evolution, trading transfer from memory against possibly
redundant computations on cached data.

Among the several options for matrix-free algorithms (see, e.g.,~\cite{Bauer2019} for an alternative), this contribution
considers the cellwise evaluation of the underlying finite-element integrals
by numerical quadrature. A general yet efficient way to matrix-free evaluation of
standard continuous and discontinuous finite element methods and tensor-product elements with tensor-product quadrature points was presented in~\cite{Anderson2021, Kronbichler2012, Kronbichler2019} for CPUs and~\cite{Anderson2021, Chalmers2023, Kolev2021, Ljungkvist2017} for GPUs. Similar developments have been made for simplicial element shapes, albeit the increased complexity due to the absence of a tensor product leads to a lower performance than for hexahedral elements~\cite{Moxey2020, Sun2020}.
Unfitted finite element methods, however, are still dominated by sparse-matrix-based implementations. The challenge lies in the loss of structure for the quadrature on intersected cells and the additional computations due to the necessary stabilization of small cut cells. The present work aims to challenge the state-of-the-art by developing techniques to efficiently implement typical terms in the weak formulation of unfitted discretizations and quantifying the respective performance.

Unfitted finite element methods are motivated by the advantage of sharp moving-interface representation and simple geometry embedding. Typically, the geometry is embedded into a Cartesian background mesh by e.g.~a level-set function, where the zero-isocontour describes the boundaries of the domain. As the background grid can be intersected arbitrarily, cells may be barely cut by the interface, leading to small volume fractions of cells that intersect with the physical domain. The small support on the shape functions leads to ill-conditioned element stiffness matrices and, therefore, to an ill-conditioning of the whole system of equations. This problem is well-known in the literature as the \textit{small cut cell problem}~\cite{Burman2010a, Prenter2018}.
A remedy for the \textit{small cut cell problem} is provided e.g.~by ghost penalty stabilization typically used in CutFEM, see~\cite{Prenter2023} for an overview of stabilization methods. Possible realizations of ghost penalty are discussed in~\cite{Guerkan2019} like the well-known face-based ghost penalty~\cite{Burman2010a} and the volume-based ghost penalty~\cite{Preuss2018}.
Dirichlet boundary conditions are usually imposed weakly in the unfitted setting. The Nitsche method~\cite{Nitsche1971} is often preferred over the penalty methods, as it adds consistent terms to the weak formulation to maintain full accuracy.
All unfitted methods have in common that integrals over cut elements need to be computed with specially adapted quadrature. There exist several methods of quadrature generation in intersected elements like octree subdivision~\cite{Duester2008}, moment fitting techniques~\cite{Sudhakar2013}, or the dimension reduction of integrals~\cite{Saye2015}.
The dimension reduction approach enables high order accuracy at the cost of unstructured quadrature. In contrast to tensor-product quadrature, for unstructured quadrature, the position of quadrature points is arbitrary inside the cell, which necessitates a separate evaluation of the shape functions for every cell and prevents the use of sum-factorization techniques~\cite{Orszag1980}. Also, the number of quadrature points per cell can vary. Therefore, an algorithm dealing with unstructured quadrature has to be more general than one for structured quadrature.

In this work, we embed a geometry into a Cartesian background mesh with tensor-product elements via a level-set function. As a remedy to the small cut cell problem, we select the volume-based ghost penalty for our study. To generate high-order quadrature rules, we employ the dimension reduction approach of~\cite{Saye2015}.
Using these methods, the main goal of this article is to bring advanced, unfitted methods together with state-of-the-art high-performance implementation techniques suited for modern hardware and massively parallel supercomputers. The implementation of the infrastructure analyzed in this work is publicly available in the widely used finite element library deal.II~\cite{Arndt2023, Arndt2021}; however, the proposed concepts are generic and would be applicable to most modern finite element realizations.

The outline of this work is as follows. Section~\ref{section:ufem} identifies integrals in the weak form to be computed with structured and unstructured quadrature, respectively. Section~\ref{section:matrixfree} details the matrix-free implementation of the presented unfitted finite element methods. The obtained accuracy is shown in Section~\ref{section:accuracy}. In Section~\ref{section:performance}, we present the computational performance on the node level for selected CPU architectures. In Section~\ref{section:application}, we analyze matrix-free operator evaluation for challenging application cases. We conclude this work with the findings of the presented work in Section~\ref{section:conclusion}.

\section{Unfitted finite element methods}\label{section:ufem}

\subsection{Discretization}

We consider a bounded domain $\Omega\subset \mathbb{R}^d$ with arbitrary shape, defined by the zero-isocontour of a level-set function.
As a model problem, consider the Poisson equation
\begin{equation}\label{eq:poisson}
  -\Laplace{u} = f\quad\mathrm{in}\;\Omega,\;\mathrm{with}\; u=g\;\mathrm{on}\;\Gamma\doteq\partial\Omega
\end{equation}
with $u$ the primary variable, $f$ the right-hand side and $g$ the Dirichlet boundary condition (DBC) on the immersed boundary $\Gamma$. We obtain the weak formulation
\begin{equation}\label{eq:weakformulation}
  a(v,u)=b(v) \qquad \text{for all test functions } v,
\end{equation}
where $a(v,u)$ is the bilinear form and $b(v)$ is the linear form. We discretize the problem with either the continuous Galerkin (CG) method in a finite-dimensional subspace of $H^1(\Omega)$ or the discontinuous Galerkin (DG) method in a finite-dimensional subspace of $L^2(\Omega)$---useful in the context of e.g.~variable coefficients or mixed finite element problems of Darcy or Stokes type~\cite{Chidyagwai2010}. In both cases, the finite-dimensional subspace is defined by tensor-product polynomial spaces on hexahedral elements in $d=3$ space dimensions.

Discretizing the problem~\eqref{eq:poisson} on an approximation of the domain $\Omega_h$ with the CG method and enforcing the DBC weakly via the Nitsche method~\cite{Nitsche1971} yields the following bilinear form:
\begin{equation}
  \begin{split}
    a_\mathrm{CG}(v,u)=&\underbrace{\intdom{\Grad{v}}{\Grad{u}}}_{\mathrm{Poisson}}
    \underbrace{-\intdomDirichlet{\partial_n v}{u}-\intdomDirichlet{v}{\partial_n u - \tau_\mathrm{D}h^{-1} u}}_{\mathrm{Nitsche}},
  \end{split}
  \label{eq:bilinear_cg}
\end{equation}
where $\partial_n = n \cdot \nabla$ denotes the normal derivative along the direction of the unit outer normal $n$ and $h$ the characteristic element scaling. The linear form reads
\begin{equation}
  b_\mathrm{CG}(v)=\underbrace{\intdom{v}{f}}_{\mathrm{Inhomogeneity}}
  \underbrace{-\intdomDirichlet{\partial_n v}{g}+\tau_\mathrm{D}h^{-1}\intdomDirichlet{v}{g}}_{\mathrm{Nitsche}}
\end{equation}
and includes the terms arising from the inhomogenity $f$ and the Nitsche terms for the DBC $g$ with the Dirichlet penalty parameter $\tau_\mathrm{D}$.

For discretizing problem~\eqref{eq:poisson} with the DG method, we use the symmetric interior penalty Galerkin (SIPG) method~\cite{Arnold2002} on interior faces and enforce the DBC weakly via the Nitsche method, leading to the bilinear form
\begin{equation}
  a_\mathrm{DG}(v,u)
  =a_\mathrm{CG}(v,u)
  \underbrace{-\intinteriorfaces{\avg{\partial_n v}}{\jump{u}}-\intinteriorfaces{\jump{v}}{\avg{\partial_n u}-\gamma h^{-1}{\jump{u}}}}_{\mathrm{SIPG}},
  \label{eq:bilinear_dg}
\end{equation}
where $\jump{\cdot}=(\cdot)^--(\cdot)^+$ denotes the jump operator, $\avg{\cdot}=\frac{1}{2}\left((\cdot)^{-}+(\cdot)^{+}\right)$ the average operator at an interior face $\Gamma_{h}^{\mathrm{int}}$, $\gamma$ the interior penalty parameter and $h$ the characteristic element scaling.
The linear form for the DG discretization
is equivalent to the CG method ($b_\mathrm{DG}(v)=b_\mathrm{CG}(v)$), apart from the different underlying function spaces.

\subsection{Ghost penalty stabilization}

The obtained weak formulation~\eqref{eq:weakformulation}--\eqref{eq:bilinear_dg} is numerically unstable as it suffers from the \textit{small cut cell problem}. The challenge is twofold: the Nitsche parameter to enforce the boundary condition is not bounded (stability), and the stiffness matrices can be arbitrarily ill-conditioned if the shape functions have only small support (conditioning). To overcome this stability and conditioning issue, stabilization approaches have been introduced in the literature. We consider the volume-based ghost penalty method~\cite{Guerkan2019, Preuss2018} throughout this work. The interested reader is referred to the supplementary materials, where we comment on the face-based ghost penalty in the context of matrix-free methods (see~SM1). Ghost penalty stabilizations add a consistent term to the bilinear form
\begin{equation}
  a_\mathrm{stab}(v,u)=a(v,u)+g_\mathrm{v}(v,u)
\end{equation}
to ensure coercivity.

For volume-based ghost penalty, we build patches $P$ over two adjacent elements that share a face and penalize the jump in the patch volume denoted by $\jump{\cdot}_P$ with
\begin{equation}
  g_\mathrm{v}(v,u)=\sum_{P\in\mathcal{P}}\frac{\tau_\mathrm{v}}{h^2}\intstabilizedpatch{\jump{v}_P}{\jump{u}_P},
  \label{eq:volume_based_gp}
\end{equation}
where $\tau_\mathrm{v}$ denotes the penalty parameter. Two patch selection strategies are common. 
The first option is to stabilize over all faces between intersected cells and physical neighbors (all neighbors strategy), the other option is to only stabilize an intersected cell across the face to the most stable neighbor (stable neighbor strategy). Stabilizing over too many faces can lead to locking effects~\cite{Burman2022} and is computationally more expensive, thus we choose the stable neighbor strategy throughout this work.

The discretized and stabilized weak formulation results in a linear system of equations of the form $\Matrix{A}\Vector{x}=\Vector{b}$. For large-scale problems like those considered in the present work and especially in three space dimensions, iterative solvers are usually the preferred choice due to the reduced complexity compared to sparse direct solvers.

\subsection{Quadrature classification}
\begin{figure}
  \centering
  \begin{tikzpicture}
    \draw[step=2.0] (0,0) grid (4,2);

    \draw (1,1) node{inside cell};

    \draw (2.8,1) node{cut cell};

    \draw (1,2.3) node{inside face};

    \draw (2.8,2.3) node{cut face};

    \draw plot[smooth, tension=1.0] coordinates {(3.1,2) (3.5, 1.3) (3.8,0.0)};

    \draw [bend left]  (3.5, 1.3) to +(1.0,0.7) node[right]{surface};

  \end{tikzpicture}
  \begin{tikzpicture}
    \draw[step=2.0] (0,0) grid (4,2);

    \draw plot[smooth, tension=1.0] coordinates {(3.1,2) (3.5, 1.3) (3.8,0.0)};

    \draw [bend left]  (3.5, 1.3) to +(1.0,0.7) node[right]{$\Gamma$};

    \foreach \px in {0.5,1.5}
    \foreach \py in {0.5,1.5}
    \filldraw[tumred] (\px, \py) circle (2pt);

    \foreach \px in {0.5,1.5,2.5,3.5}
    \foreach \py in {0.0,2.0}
    \filldraw[tumorange] (\px, \py) circle (2pt);
    \foreach \px in {0.0,2.0}
    \foreach \py in {0.5,1.5}
    \filldraw[tumorange] (\px, \py) circle (2pt);

    \foreach \px in {2.4,3.4}
    \draw[tumred] (\px, 0.5) circle (2pt);
    \foreach \px in {2.3,3.1}
    \draw[tumred] (\px, 1.5) circle (2pt);

    \draw[tumblue] (3.71, 0.5) circle (2pt);
    \draw[tumblue] (3.42, 1.5) circle (2pt);

    \foreach \px in {2.3,3.3}
    \draw[tumorange] (\px, 0.0) circle (2pt);
    \foreach \px in {2.2,2.9}
    \draw[tumorange] (\px, 2.0) circle (2pt);

  \end{tikzpicture}
  \caption{Classification of weak formulation terms in quadrature on respective cell/face type: volume quadrature (red), face quadrature (orange) and surface quadrature (blue); structured (tensor-product) quadrature (dots) and unstructured quadrature (circles).}\label{fig:quadratures}
\end{figure}
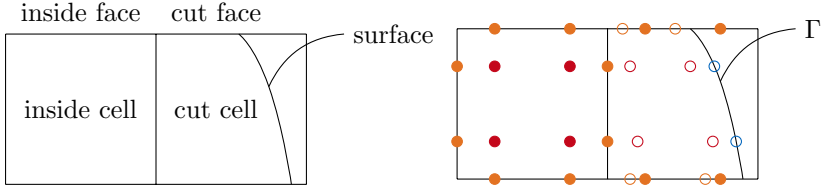

 We distinguish between two types of quadrature in this work: tensor-product quadrature, called structured quadrature throughout this work, and non-tensor-product quadrature, called unstructured quadrature.
When evaluating the terms arising in the weak form, the following general classification with respect to the structure of the quadrature can be made:
\begin{itemize}
  \item Poisson term, inhomogenity:
  \begin{itemize}
    \item inside cells $\rightarrow$ structured cell quadrature
    \item cut cells $\rightarrow$ unstructured cell quadrature
  \end{itemize}
  \item Nitsche term: unstructured surface quadrature (only on cut cells)
  \item SIPG term: 
  \begin{itemize}
    \item inside face $\rightarrow$ structured face quadrature
    \item cut face $\rightarrow$ unstructured face quadrature
  \end{itemize}
  \item volume-based ghost penalty: structured cell quadrature
\end{itemize}
This classification allows the use of optimal complexity algorithms for the respective term. The different kinds of quadrature are visualized in Figure~\ref{fig:quadratures}.

\section{Matrix-free implementation}\label{section:matrixfree}

\subsection{Matrix-free vs.\ matrix-based method}\label{subsection:matrixfreeandbased}

One main component of iterative linear solvers is the operator action
\begin{equation}\label{eq:operatorevaluation}
  \Vector{y}=A(\Vector{u}),
\end{equation}
which is a more abstract notation of a matrix-vector product. Here, we denote by $\Vector{u}$ the source vector and by $\Vector{y}$ the result (destination) vector of the operator.
In the proposed matrix-free algorithms, the operator evaluation is defined by a loop over
cells and faces. Even though different ways to schedule the face
integrals in a DG method are possible~\cite{Kronbichler2019}, this work only
considers the case that runs separate loops over cells and faces with the goal to
amortize as much of the arithmetic work as possible when evaluating the numerical fluxes appearing in the
face terms from the two adjacent sides. The operator
evaluation is thus defined as
\begin{equation}\label{eq:matrixfreeopcompact}
  A(\Vector{u})
  =\sum_{e=1}^{n_\mathrm{ele}}A_{e}(\Vector{u})+\sum_{f=1}^{n_\mathrm{faces}}A_{f}(\Vector{u})
  =\left(\sum_{e=1}^{n_\mathrm{ele}}D_e I_e B_e E_e R_e  + \sum_{f=1}^{n_\mathrm{faces}}D_f I_f B_f E_f R_f\right) \Vector{u},
\end{equation}
where the following ingredients are executed in sequence for elements and faces, respectively:
\begin{itemize}
  \item read operator $R_*$: Extract local DoFs from the global vector.
  \item evaluate operator $E_*$: Evaluate interpolation of local solution and/or its derivative in reference coordinates at quadrature points.
  \item differential operator $B_*$: At each quadrature point, apply metric terms and evaluate quantities of weak form.
  \item integrate operator $I_*$: Multiply with all test functions relevant to cell/face and accumulate integral contributions across all quadrature points. For the self-adjoint Poisson operator, $I_*$ is the transpose of $E_*$.
  \item distribute operator $D_*$: Add local contributions to global vector/matrix. It is the transposed operator to read for the Poisson problem.
\end{itemize}
The above five cell-/facewise operators are the basis for both matrix-free operator evaluation and for assembling the local matrices $\Matrix{A}_e = I_e B_e E_e$ into a global sparse matrix. We reference them throughout this work where appropriate.

When denoting by $\Matrix{A}$ the matrix representing the operator $A(\Vector{u})$, the evaluation~\eqref{eq:matrixfreeopcompact} is equivalent to the matrix-vector product
\begin{equation}\label{eq:matrixbasedop}
\Vector{y}=\Matrix{A}\Vector{u}.
\end{equation}
For a software implementation that does not build matrix representations of the three operators $E_e, B_e, I_e$ but only applies their respective action on vectors, the local matrices $\Matrix{A}_e$ can also be computed column-by-column by applying $I_e B_e E_e$  on element-wise unit vectors.

To the authors' knowledge, all current implementations of unfitted FEM are
based on sparse matrices, which have the drawback of an exploding complexity
per degree of freedom (DoF) as the polynomial order increases due to the dense coupling between all DoFs on (adjacent) elements. While matrix-free methods also
exhibit the coupling between the unknowns, they obviate accessing
the full $n_\mathrm{DoF,e}^2$ matrix entries per cell by amortizing the data
access for all $n_\mathrm{DoF,e}$ DoFs with repeated computations on
cached data, which fits with the high arithmetic capabilities
of contemporary hardware. Hence, computing the cell (and face) integrals for every
operator application~\eqref{eq:matrixfreeopcompact} can be attractive
for high-order elements, despite redundant arithmetic work. The
subsequent sections will identify algorithms with the best balance between
memory access and arithmetic cost, thus minimizing the time for
operator evaluation.

\subsection{Evaluation of different quadrature}
To use optimal complexity algorithms, we have classified the terms arising from the weak formulation with respect to the underlying quadrature point distribution into the \textit{structured} and \textit{unstructured} cases. To reflect this choice,
we can split up the operator evaluation~\eqref{eq:operatorevaluation} into
\begin{equation}\label{eq:additive_decomp}
  \Vector{y}=\Vector{y}_\mathrm{unstr}+\Vector{y}_\mathrm{str},
\end{equation}
involving contributions from structured quadrature $\Vector{y}_\mathrm{str}$ and unstructured quadrature $\Vector{y}_\mathrm{unstr}$. This leads to the following options for implementing the operator evaluation:
\begin{center}
\begin{forest}
  [operator evaluation
  [matrix-based]
  [mixed
  ($\Vector{y}_\mathrm{unstr}$ matrix-based + $\Vector{y}_\mathrm{str}$ matrix-free)
  ]
  [matrix-free]
  ]
\end{forest}
\end{center}
For \textit{structured} quadrature, it has been shown that matrix-free methods are faster than sparse matrices~\cite{Kolev2021, Kronbichler2019, Kronbichler2018}. Therefore, a matrix-free approach is preferable for all the structured terms according to the additive decomposition~\eqref{eq:additive_decomp}. We have two options for the terms with \textit{unstructured} quadrature: precompute or recompute. Precomputing means we evaluate the unstructured integrals in a setup step and store the contributions in either a sparse matrix or the element matrices, see also~\cite{Cantwell2011} for the general concepts. Recomputing means we also evaluate the underlying integrals in every operator evaluation. In the following, we present efficient algorithms for these two steps.

\begin{remark}
  We use Einstein notation where appropriate for compact notation. We show formulas either dimension independent or for $d=3$; for clarity, figures are for $d=2$.
\end{remark}

As a starting point for the optimized algorithms for the operators $E_*$ and $I_*$, consider the basic evaluation $E_*$ of the finite element interpolation $u$ at a quadrature point, which we call $\hat{u}$,
\begin{equation}
  \hat{u}_q 
  =\sum_{l=1}^{n_\mathrm{DoF}^\mathrm{dim}}\varphi_l(\Vector{x}_q)u_l
  =S_{ql}u_l\quad \mathrm{with}\;q=0,...,n_\mathrm{q}^\mathrm{dim},
  \label{eq:naive_evaluation}
\end{equation}
where $n_\mathrm{DoF}^\mathrm{dim}$ is the number of DoFs on a cell, $n_\mathrm{q}^\mathrm{dim}$ the number of quadrature points, and $S$ denotes the interpolation matrix of shape functions evaluated at the quadrature points. The arithmetic complexity of naive interpolation assuming $n_\mathrm{DoF}^\mathrm{dim}=n_\mathrm{q}^\mathrm{dim}=k^d$ is $\mathcal{O}(k^{2d})$ and the memory complexity for the interpolation matrix $S$ is $\mathcal{O}(k^{2d})$. In this work, we choose the 1D number of quadrature points $n_\mathrm{q}^{1\mathrm{D}}$ as $k=p+1$ where $p$ is the polynomial degree of the shape functions. For a tensor-product element, the number of DoFs per direction $n_\mathrm{DoF}^{1\mathrm{D}}$ is also equal to $k$.

\subsubsection*{Structured quadrature evaluation}

For inside cells, the quadrature involves points with a tensor-product structure identical for every cell. We use sum factorization~\cite{Orszag1980} with optimized implementations as developed in~\cite{Kronbichler2012, Kronbichler2019} for this cell type. The implementation uses explicit data-level parallelism through the single-instruction/multiple-data (SIMD) paradigm by processing the integrals of multiple cells in different lanes of a single instruction. This forms batches of cells, with a batch consisting of as many cells as the hardware-dependent width of SIMD vectors.

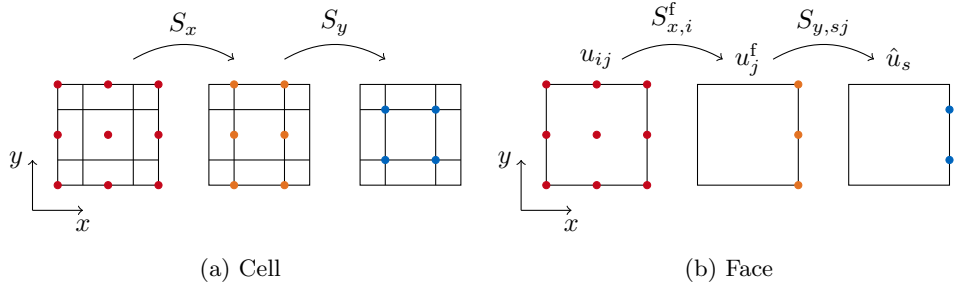
\begin{figure}
  \centering
  \begin{subfigure}[b]{0.49\textwidth}
    \begin{tikzpicture}[scale=0.67]
      \draw[->] (-0.5,-0.5) -- +(1.0,0) node[below]{$x$};
      \draw[->] (-0.5,-0.5) -- +(0,1.0) node[left]{$y$};
  
      \def \elelength {2};
      \def \stride {1};
  
      \foreach \x in {0,1,2}
      {
        \draw (\stride*\x+\elelength*\x, 0) rectangle +(\elelength, \elelength);
      }
      \foreach \y in {0,1}
        \draw (0, 0.25*\elelength+0.5*\y*\elelength) -- +(\elelength, 0);
      \foreach \x in {0,1}
        \draw (0.25*\elelength+0.5*\x*\elelength, 0) -- +(0, \elelength);
      \foreach \px in {0,1,2}
      \foreach \py in {0,1,2}
      \filldraw[tumred] (\px, \py) circle (2pt);
  
      \draw [bend left, ->]  (0.75*\elelength,\elelength + 0.5) to  node[above]{$S_x$} +(\stride+0.5*\elelength,0);
  
      \foreach \y in {0,1}
        \draw (\elelength+\stride, 0.25*\elelength+0.5*\y*\elelength) -- +(\elelength, 0);
      \foreach \x in {0,1}
        \draw (1.25*\elelength+0.5*\x*\elelength+\stride, 0) -- +(0, \elelength);
      \foreach \x in {0,1}
      \foreach \py in {0,1,2}
        \filldraw[tumorange] (1.25*\elelength+0.5*\x*\elelength+\stride, \py) circle (2pt);
  
      \draw [bend left, ->]  (1.75*\elelength + \stride,\elelength + 0.5) to  node[above]{$S_y$} +(\stride+0.5*\elelength,0);
  
      \foreach \y in {0,1}
        \draw (2*\elelength+2*\stride, 0.25*\elelength+0.5*\y*\elelength) -- +(\elelength, 0);
      \foreach \x in {0,1}
        \draw (2.25*\elelength+0.5*\x*\elelength+2*\stride, 0) -- +(0, \elelength);
      \foreach \x in {0,1}
      \foreach \y in {0,1}
        \filldraw[tumblue] (2.25*\elelength+0.5*\x*\elelength+2*\stride, 0.25*\elelength+0.5*\y*\elelength) circle (2pt);
     
    \end{tikzpicture}
    \caption{Cell}\label{subfig:structured_cell}
  \end{subfigure}
  \begin{subfigure}[b]{0.49\textwidth}
    \begin{tikzpicture}[scale=0.67]
      \draw[->] (-0.5,-0.5) -- +(1.0,0) node[below]{$x$};
      \draw[->] (-0.5,-0.5) -- +(0,1.0) node[left]{$y$};
  
      \def \elelength {2};
      \def \stride {1};
  
      \foreach \x in {0,1,2}
      {
        \draw (\stride*\x+\elelength*\x, 0) rectangle +(\elelength, \elelength);
      }
  
      \foreach \px in {0,1,2}
      \foreach \py in {0,1,2}
      \filldraw[tumred] (\px, \py) circle (2pt);
  
      \draw(0.5*\elelength,\elelength + 0.5) node{$u_{ij}$};
  
      \draw [bend left, ->]  (0.75*\elelength,\elelength + 0.5) to  node[above]{$S^\mathrm{f}_{x,i}$} +(\stride+0.5*\elelength,0);
  
      \draw (1.5*\elelength + \stride,\elelength + 0.5) node{$u^\mathrm{f}_{j}$};
  
      \foreach \py in {0,1,2}
        \filldraw[tumorange] (2*\elelength+\stride, \py) circle (2pt);
  
      \draw [bend left, ->]  (1.75*\elelength + \stride,\elelength + 0.5) to  node[above]{$S_{y,sj}$} +(\stride+0.5*\elelength,0);
  
      \draw (2.5*\elelength + 2*\stride,\elelength + 0.5) node{$\hat{u}_{s}$};
  
      \foreach \y in {0,1}
        \filldraw[tumblue] (3*\elelength + 2*\stride, 0.25*\elelength+0.5*\y*\elelength) circle (2pt);
     
    \end{tikzpicture}
    \caption{Face}\label{subfig:structured_face}
  \end{subfigure}
  \caption{Structured cell and face tensor-product quadrature evaluation using sum-factorization for a 2D element with $p=2$. \subref{subfig:structured_cell} Interpolation from DoF values into tensor-product quadrature points using the tensor-product structure of the shape functions. \subref{subfig:structured_face} Interpolation from cell DoF values into face DoF values, then in-face interpolation into quadrature points.}\label{fig:structured}
\end{figure}

To evaluate~\eqref{eq:naive_evaluation} at all quadrature points with tensor indices $(r,s,t)$, we use the tensor-product structure of shape functions and quadrature points to factor out the evaluation in a single direction at a time,
\begin{equation}
    \hat{u}_q 
    = \sum_{l=1}^{n_\mathrm{DoF}^\mathrm{3D}}\varphi_l(\Vector{x}_q)u_l
    =\hat{u}_{rst} 
    =\sum_{i=1}^{n_\mathrm{DoF}^\mathrm{1D}}\varphi_i^\mathrm{1D}(x_{r}^\mathrm{1D})\sum_{j=1}^{n_\mathrm{DoF}^\mathrm{1D}}\varphi_j^\mathrm{1D}(y_{s}^\mathrm{1D})\sum_{k=1}^{n_\mathrm{DoF}^\mathrm{1D}}\varphi_k^\mathrm{1D}(z_{t}^\mathrm{1D})u_{ijk},
  \end{equation}
with $l = k \left({n_\mathrm{DoF}^\mathrm{1D}}\right)^2 + j\, n_\mathrm{DoF}^\mathrm{1D} + i$ and
$q = t \left({n_\mathrm{q}^\mathrm{1D}}\right)^2 + s\, n_\mathrm{q}^\mathrm{1D} + r$.
In the optimized evaluation with sum factorization, each sum is evaluated on the entirety of the index space, e.g., for the sum over $i$, all indices $j,k$ are evaluated for all indices $r$, before continuing with the summation over $j$ and $k$, respectively. We write
\begin{equation}
  \hat{u}_{rst} = S_{z,tk} S_{y,sj} S_{x,ri} u_{ijk},
  \label{eq:sum_factorization_cell}
\end{equation}
where $S_{*,\star}$ denotes the matrix of 1D shape functions evaluated at the 1D position of quadrature points in direction $*$ (dimension $n_\mathrm{q}^\mathrm{1D}\times n_\mathrm{DoF}^\mathrm{1D}$), doing one contraction at a time (see Figure~\ref{subfig:structured_cell}). This reduces the complexity from $\mathcal{O}(k^{2d})$ (of the naive implementation) to $\mathcal{O}(dk^{d+1})$. The integration $I_e$ as the transposed operation $\hat{u}_{rst}\rightarrow u_{ijk}$ is done analogously, and similarly for the gradients~\cite{Kronbichler2019}.

For SIPG terms on inside faces (or face-based ghost penalty terms) with structured tensor-product quadrature, we can use the efficient techniques presented in~\cite{Kronbichler2019}. The procedure is similar to the cell case. First, values (see Equation~\eqref{eq:sum_factorization_face}) and the normal derivative are projected onto the faces DoFs $u^\mathrm{f}$, and then inside the face, sum-factorization is used in dimension $d-1$, see Figure~\ref{subfig:structured_face}. We define the evaluation of quadrature values $\hat{u}$ on faces $E_f$ for the example of a face aligned with the $x$-direction as
\begin{equation}
  \hat{u}_{st} = S_{z,tk} S_{y,sj} \underbrace{S_{x,i}^\mathrm{f} u_{ijk}}_{u^\mathrm{f}_{jk}},
  \label{eq:sum_factorization_face}
\end{equation}
where $u^\mathrm{f}$ are the DoF values on the face. Note that for nodal Lagrangian elements defined in the Gauss--Lobatto nodes, no interpolation is necessary for the values, as the face DoF values $u^\mathrm{f}$ can be read from the cell DoF values $u$ with the $\delta_{ij}$ property. For the gradients in the normal direction to the face, however, an interpolation into the face DoFs is necessary.
Again, integration $I_f$ is implemented analogously. The arithmetic complexity for structured face evaluation is $\mathcal{O}(k^d)$ per face.

\subsubsection*{Unstructured quadrature evaluation}

\begin{figure}
  \centering
  \begin{subfigure}[b]{0.49\textwidth}
    \begin{tikzpicture}[scale=0.67]
      \draw[->] (-0.5,-0.5) -- +(1.0,0) node[below]{$x$};
      \draw[->] (-0.5,-0.5) -- +(0,1.0) node[left]{$y$};
  
      \def \elelength {2};
      \def \stride {1};
  
      \foreach \x in {0,1,2}
      {
        \draw (\stride*\x+\elelength*\x, 0) rectangle +(\elelength, \elelength);
      }
      \draw (0, 0.75*\elelength) -- +(\elelength, 0);
      \draw (0.75*\elelength, 0) -- +(0, \elelength);
      \foreach \px in {0,1,2}
      \foreach \py in {0,1,2}
      \filldraw[tumred] (\px, \py) circle (2pt);
  
      \draw [bend left, ->]  (0.75*\elelength,\elelength + 0.5) to  node[above]{$S^q_x$} +(\stride+0.5*\elelength,0);
  
      \draw (\elelength+\stride, 0.75*\elelength) -- +(\elelength, 0);
      \draw (1.75*\elelength+\stride, 0) -- +(0, \elelength);
      \foreach \py in {0,1,2}
      \filldraw[tumorange] (1.75*\elelength+\stride, \py) circle (2pt);
  
      \draw [bend left, ->]  (1.75*\elelength + \stride,\elelength + 0.5) to  node[above]{$S^q_y$} +(\stride+0.5*\elelength,0);
  
      \draw (2*\elelength+2*\stride, 0.75*\elelength) -- +(\elelength, 0);
      \draw (2.75*\elelength+2*\stride, 0) -- +(0, \elelength);
      \filldraw[tumblue] (2.75*\elelength+2*\stride, 0.75*\elelength) circle (2pt);
     
    \end{tikzpicture}
    \caption{Cell}\label{subfig:unstructured_cell}
  \end{subfigure}
  \begin{subfigure}[b]{0.49\textwidth}
    \begin{tikzpicture}[scale=0.67]
      \draw[->] (-0.5,-0.5) -- +(1.0,0) node[below]{$x$};
      \draw[->] (-0.5,-0.5) -- +(0,1.0) node[left]{$y$};
  
      \def \elelength {2};
      \def \stride {1};
  
      \foreach \x in {0,1,2}
      {
        \draw (\stride*\x+\elelength*\x, 0) rectangle +(\elelength, \elelength);
      }
  
      \foreach \px in {0,1,2}
      \foreach \py in {0,1,2}
      \filldraw[tumred] (\px, \py) circle (2pt);
  
      \draw(0.5*\elelength,\elelength + 0.5) node{$u_{ij}$};
  
      \draw [bend left, ->]  (0.75*\elelength,\elelength + 0.5) to  node[above]{$S^\mathrm{f}_{x,i}$} +(\stride+0.5*\elelength,0);
  
      \draw (1.5*\elelength + \stride,\elelength + 0.5) node{$u^\mathrm{f}_{j}$};
  
      \foreach \py in {0,1,2}
        \filldraw[tumorange] (2*\elelength+\stride, \py) circle (2pt);
  
      \draw [bend left, ->]  (1.75*\elelength + \stride,\elelength + 0.5) to  node[above]{$S^q_{y,j}$} +(\stride+0.5*\elelength,0);
  
      \draw (2.5*\elelength + 2*\stride,\elelength + 0.5) node{$\hat{u}_{q}$};
  
      \foreach \y in {0}
        \filldraw[tumblue] (3*\elelength + 2*\stride, 0.75*\elelength+0.5*\y*\elelength) circle (2pt);
     
    \end{tikzpicture}
    \caption{Face}\label{subfig:unstructured_face}
  \end{subfigure}
  \caption{Unstructured cell and face quadrature evaluation using the tensor-product structure of the shape functions for a 2D element with $p=2$. \subref{subfig:unstructured_cell} Interpolation from DoF values into an arbitrary quadrature point using 1D interpolations. \subref{subfig:unstructured_face} Interpolation from the cell DoFs into the face DoFs, then in-face interpolation an arbitrary quadrature point on the face}\label{fig:unstructured}
\end{figure}
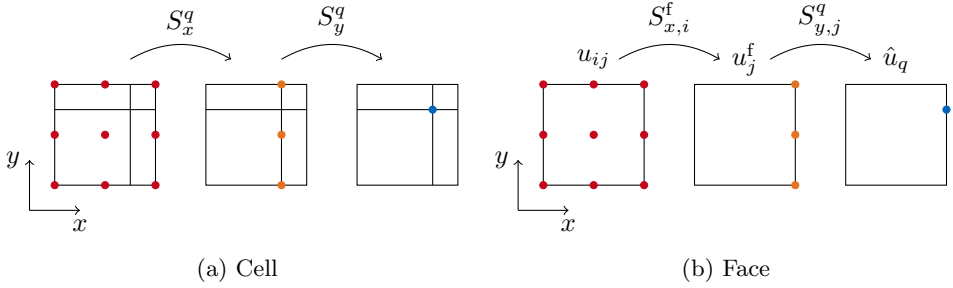

For the matrix-free evaluation with unstructured quadrature points on cells (Poisson and Nitsche terms) and faces (SIPG term) that arise in unfitted finite element type discretizations, we implement an efficient algorithm that makes use of the tensor-product structure of the underlying finite element shape functions on the hexahedral element, but working on a single point at a time,
\begin{equation}
    \hat{u}_q = \sum_{l=1}^{n_\mathrm{DoF}^\mathrm{3D}}\varphi_l(\Vector{x}_q)u_l=\sum_{i=1}^{n_\mathrm{DoF}^\mathrm{1D}}\varphi_i^\mathrm{1D}(x_q)\sum_{j=1}^{n_\mathrm{DoF}^\mathrm{1D}}\varphi_j^\mathrm{1D}(y_q) \sum_{k=1}^{n_\mathrm{DoF}^\mathrm{1D}}\varphi_k^\mathrm{1D}(z_q) u_{ijk}
\end{equation}
 with $l = k \left({n_\mathrm{DoF}^\mathrm{1D}}\right)^2 + j\, n_\mathrm{DoF}^\mathrm{1D} + i$. Here, we split the sum over all basis functions into $d$ one-dimensional sums exploiting the tensor-product structure of the basis functions. In compact notation, we write
\begin{equation}
  \hat{u}_q = S^q_{z,k} S^q_{y,j} S^q_{x,i} u_{ijk},
  \label{eq:unstructured_cell}
\end{equation}
where $S^q_{*,\star}$ is the vector of evaluated 1D basis functions in the respective direction $*$ for quadrature point $q$, see Figure~\ref{subfig:unstructured_cell}. Equation~\eqref{eq:unstructured_cell} implies an arithmetic complexity of evaluation of unstructured quadrature for all $\mathcal{O}(k^d)$ quadrature points of $\mathcal{O}(k^{2d})$, which is relatively high compared to the structured case with sum factorization.

One key ingredient with tensor-product shape functions is the lower arithmetic complexity to evaluate the actual polynomial basis by 1D evaluations, leading to a $\mathcal{O}(2dk(k-1))$ complexity in one-dimensional DoFs $k$ of $d$ dimensions, leading to significant savings for high $k$ as opposed to evaluating generic $d$-dimensional polynomials.
We precompute the 1D polynomial evaluation in an initialization step and reuse the data during evaluation $E$ and integration $I$. The reduced data complexity of 1D vectors $\mathcal{O}(dk)$ vs.~the expanded $d$-tensor data $\mathcal{O}(k^d)$ used by the naive interpolation~\eqref{eq:naive_evaluation} is crucial to keep the data in the low-level caches for fast access.

DG methods need to compute the numerical fluxes only on the physically relevant part of faces between cut cells. The unstructured quadrature points on the face necessitate specialized evaluation techniques as well. The algorithm utilizes the knowledge that all quadrature points are located on the face to reduce the dimension in which the unstructured quadrature algorithm needs to be run.

The first step is to interpolate values and normal derivatives on the face DoFs $u^\mathrm{f}$, keeping the tensor product within the face. This step is equivalent to the structured quadrature case (see Figure~\ref{subfig:unstructured_face} and Equation~\ref{eq:sum_factorization_face}). For each quadrature point, a $d-1$ dimensional interpolation in-face exemplarily shown for a $x$-direction face
\begin{equation}
  \hat{u}_q = S^q_{z,k} S^q_{y,j} \underbrace{S_{x,i}^\mathrm{f} u_{ijk}}_{u^\mathrm{f}_{jk}}
  \label{eq:unstructured_face}
\end{equation}
is performed. The arithmetic complexity of unstructured face evaluation at all quadrature points of the face is  $\mathcal{O}(k^{2(d-1)})$.

Following the rationale of previous work that found that automatic vectorization by the compiler does not adequately use the compute capabilities~\cite{Kronbichler2012, Kronbichler2019}, we explicitly implement the SIMD vectorization by intrinsics in an outer-loop fashion, combining the operations of several quadrature points of an entity into the SIMD lanes. The metric terms derived from the map between reference and real coordinates, i.e.,~the inverse Jacobians and the Jacobian determinant multiplied by the quadrature weight, are precomputed to save computational effort and loaded in a vectorized (packed) fashion from aligned memory.

\begin{remark}
  The class performing the presented structured quadrature algorithm in deal.II, which is used for the numerical experiments below, is called \texttt{FEEvaluation}, and the class for unstructured quadrature evaluation is called \texttt{FEPointEvaluation}. The latter is the major contribution of this work. We refer the reader to the code and documentation for further information and interpretation of the above formulas.
\end{remark}

\begin{remark}
  The intrinsics-based vectorization over points is accomplished by using an abstraction class called \texttt{VectorizedArray} \cite{Kronbichler2012} (a class similar to \texttt{std::simd}) as the underlying arithmetic data type of a point. This class overloads the basic arithmetic operations, allowing for templated code with a generic arithmetic number type that, under the hood, leverages the SIMD capabilities of current computer hardware.
\end{remark}

\subsection{Evaluation of stabilization terms}\label{subsec:eval_stab}

In addition to the physical terms in the weak formulation, the integrals for stabilization terms need to be evaluated in a matrix-free realization.

For the volume-based ghost penalty~\eqref{eq:volume_based_gp}, we split the cell patch integrals over face-neighbors into two cell integrals. As the cells are face-neighbors, we implement the term in the face loop. The data access is identical to standard face integrals for Lagrangian finite elements. In the DG case, we can reuse the data read for the flux calculation and write the result only once.

The patch volume can be split into the two element volumes $P=E_1\cup E_2$. We can then define the volume jump for each element, e.g.~for $u$ to
\begin{equation}
  \jump{u}_{E_1}=\jump{u_1-\mathcal{E}u_2}_{E_1}\quad\mathrm{and}\quad\jump{u}_{E_2}=\jump{\mathcal{E}u_1-u_2}_{E_2}
\end{equation}
where $\mathcal{E}$ defines the extrapolation operator of the solution on one element to the other element and, therefore, the whole patch. For Cartesian background grids like the ones used in the present work, we implement the extension operator $\mathcal{E}$ by shifting the quadrature points of the reference element into the neighboring element. Therefore, the quadrature rule preserves its tensor-product structure, enabling the use of the efficient evaluation routines presented in Section~\ref{section:matrixfree} with the only difference that for the extrapolation, anisotropic shape values at quadrature points (e.g. $S_x\neq S_y=S_z$ for extrapolation in $x$-direction) have to be used instead of the isotropic ones, cf.~Equation~\eqref{eq:sum_factorization_cell}.

\subsection{Combined matrix-free algorithm design}

The cell and face loops for operator evaluation~\eqref{eq:operatorevaluation} are structured as shown in Algorithms~\ref{alg:combined_cell} and~\ref{alg:combined_face}. The essential design choice is to keep as large a portion of the algorithm as possible vectorized across cells/faces because overhead can be reduced by amortizing the work of setting pointers, running loops and looking up data across several entities.

For cells (see Algorithm~\ref{alg:combined_cell}), we use a categorization to allow only cells of the same classification in a SIMD cell batch, implying there are either inside or intersected cells in one batch. The structured algorithm infrastructure handles reading from and writing to global vectors. The unstructured algorithm directly writes into the data buffer for local DoF values used by the structured algorithm infrastructure to avoid additional data copies or indirect addressing.
\begin{algorithm}[t]
  \caption{Combined structured/unstructured cell evaluation with vectorization over cells in the structured algorithm and vectorization over quadrature points in the unstructured algorithm}\label{alg:combined_cell}
  \begin{algorithmic}
    \For{cell batch = 0, ..., n cell batches}
    \State{read cell DoFs}\Comment{structured algorithm infrastructure}
      \If{inside}
        \State{do Poisson term}\Comment{structured algorithm}
      \Else
        \If{intersected}
        \For{lane = 0,...,$n_\mathrm{lanes}$}
          \State{do Poisson + Nitsche term}\Comment{unstructured algorithm}
        \EndFor
        \EndIf
      \EndIf
      \State{distribute cell DoFs}\Comment{structured algorithm infrastructure}
    \EndFor
  \end{algorithmic}
\end{algorithm}

Faces are grouped together to batches with $n_\mathrm{lanes}$ faces per batch, similarly to the cells. A batch involves faces with the same local face numbers (describing the position in the reference cell) on the respective sides, such that the face-normal interpolation is identical for the whole batch.
For faces (see Algorithm~\ref{alg:combined_face}), arbitrary combinations of inside and intersected faces are allowed in a batch. Access to the global vectors as well as the face-normal interpolation is identical for structured and unstructured evaluation because it only depends on the DoFs and not the quadrature, so the structured quadrature algorithm infrastructure can be used. For the $(d-1)$-dimensional in-face operations different paths have to be selected at run time. For the structured algorithm, lanes without contributions from structured quadrature are masked out by multiplying the respective lane values with zero. The unstructured algorithm works by accessing only specific lanes in the local buffers for the face's DoF values used.

\begin{algorithm}[t]
  \caption{Combined structured/unstructured face evaluation (DG + ghost penalty) with vectorization over faces in the structured algorithm and vectorization over quadrature points in the unstructured algorithm}\label{alg:combined_face}
  \begin{algorithmic}
    \For{face batch = 0, ..., n face batches}
    \State{check for ghost penalty, inside flux and cut flux on lanes}
    \State{read cell DoF values}\Comment{structured algorithm infrastructure}
    \State{interpolate into face DoF values}\Comment{structured algorithm, face-normal}
    \If{batch has ghost penalty or SIPG flux on inside face}
    \State{do ghost penalty term and SIPG flux}\Comment{structured algorithm, in-face}
    \EndIf
    \If{batch has SIPG flux on cut face}
    \For{lane=0,...,$n_\mathrm{lanes}$}
    \If{lane has SIPG flux on cut face}
    \State{do SIPG flux} \Comment{unstructured algorithm, in-face}
    \EndIf
    \EndFor
    \EndIf
    \State{collect from face DoF values}\Comment{structured algorithm, face-normal}
    \State{distribute cell DoF values}\Comment{structured algorithm infrastructure}
    \EndFor
  \end{algorithmic}
\end{algorithm}

\subsection{Load-balancing in parallel computations}
The necessity for effective load balancing in parallel computations arises from variations in the evaluation and integration complexities inherent to interior and intersected cells. Following the categorization of cells into distinct groups, a repartitioning of the triangulation is conducted. Achieving a desirable load balance involves the assignment of weights to individual cell categories, reflecting the relative computational costs associated with each category. By default, unless explicitly specified otherwise, we have employed a weight of 1 for interior cells, 10 for intersected cells, and 0 for exterior cells. We intentionally select a high weight for intersected cells because the number of quadrature points can vary for the unstructured quadrature on the cell and faces depending on the cut configuration. Assuming that the bulk of the elements is located inside and only a small amount of cells intersect the interface, achieving a well-balanced distribution of the bulk across multiple processors becomes crucial. Consequently, maintaining high utilization is prioritized, even if processors handling intersected elements complete their tasks ahead of others. For the matrix-based operator evaluation, the weights are chosen as 1 for interior cells, 1 for intersected cells, and 0 for exterior cells.

\subsection{Iterative solution}

We use a preconditioned conjugate gradient iteration to solve the resulting linear system. From a performance perspective, two building blocks are most relevant: the operator evaluation (matrix-vector product) and the preconditioner. This work focuses on the first part, optimizing the operator evaluation.

For the preconditioner, we use state-of-the-art techniques, albeit further improvements are planned in subsequent work. For high-order discontinuous Galerkin finite elements, we use an optimal complexity hybrid polynomial multigrid preconditioner from the open-source project \texttt{ExaDG}~\cite{Arndt2020}. This multigrid preconditioner transfers from the DG space to a continuous finite element space ($c$-transfer) down to linear polynomials ($p$-transfer)~\cite{Fehn2020, Munch2023}, see also the general solver comparison for matrix-free and matrix-based algorithms in~\cite{Kronbichler2018} and references therein. On the coarse level, we use a conjugate gradient solver preconditioned by an algebraic multigrid preconditioner for which we use the \texttt{Trilinos ML} package~\cite{Gee2006} with a Chebyshev smoother.
To overcome the ill-conditioning still present in the linear system of equations after stabilization (which is partly caused by ill-conditioned local polynomial bases when extended into neighbors~\cite{Schoeder2020}), we employ a Chebyshev smoother with a computationally expensive additive Schwarz inner preconditioner~\cite{Prenter2019} on the multigrid levels. Therein, the inverses of the diagonal blocks are applied through matrix-based Cholesky decompositions utilizing \texttt{LAPACK}~\cite{Anderson1999}. For continuous elements inside the computational domain with no contributions from cut elements or stabilization, we exploit the tensor-product structure of shape functions and quadrature points to apply the fast diagonalization method~\cite{Lottes2005, Lynch1964, Munch2023Cache} to maintain optimal computational complexity.

\section{Accuracy}\label{section:accuracy}

High-order methods are motivated by their improved convergence, allowing to reach a higher accuracy with a given budget of degrees of freedom or, more frequently, to reach a certain accuracy with fewer degrees of freedom, see, e.g.,~\cite{Kronbichler2018}. The ultimate goal is to minimize the computational time for a prescribed accuracy. We demonstrate optimal convergence orders for unfitted FEM with a sphere benchmark.

We define the physical domain $\Omega$ as the unit ball with radius $R=1$ around the origin. As the computational domain, we use a Cartesian background mesh that is unfitted but larger than the physical domain $\Vector{x}\in[-1.035,+1.035]^d$.
We choose the right-hand side term as
\begin{equation}
  f(r)=-\omega^2(\cos(\omega r) + (d - 1) \mathrm{sinc}(\omega r)),
\end{equation}
which yields the analytical solution
\begin{equation}
  u_\mathrm{analytic}(r)=\cos(\omega R)-\cos(\omega r),
\end{equation}
where $r^2=\sum_{i=1}^{d}x_i^2$ and $R$ is the radius of the sphere, where a zero Dirichlet boundary condition ($g=0$) is enforced.

We use the penalty parameter $\tau_\mathrm{v}=1$ in the volume ghost penalty method. The optimal choice of $\tau_\mathrm{v}$ is still active research.

\pgfplotsset{p1/.style={tumorange,mark=o, semithick}}
\pgfplotsset{p2/.style={tumgreen,mark=otimes, semithick}}
\pgfplotsset{p3/.style={tumred,mark=*, semithick}}
\pgfplotsset{p4/.style={diag_violet,mark=oplus, semithick}}
\begin{figure}
  \begin{subfigure}[b]{0.49\textwidth}
    \begin{tikzpicture}
      \begin{axis}[
      x=1.5cm,
      y=0.2cm,
      axis x line=bottom,
      axis y line=left,
      xmode=log,
      ymode=log,
      xlabel=h/L,ylabel= relative $L^2$ error,
      grid=major,
      legend pos=outer north east,
      ]
      \addplot [domain=1.0e-2:2.0e-1, samples=100, dashed
      ]{x^2};
      \addplot [domain=1.0e-2:2.0e-1, samples=100, dashed
      ]{0.3*x^3};
      \addplot [domain=2.0e-2:2.0e-1, samples=100, dashed
      ]{0.1*x^4};
      \addplot [domain=4.0e-2:2.0e-1, samples=100, dashed
      ]{0.1*x^5};
      \addplot [p1
      ] coordinates {
          (1.725000000e-01, 3.283013484e-02)
          (8.625000000e-02, 7.893285973e-03)
          (4.312500000e-02, 1.898945659e-03)
          (2.156250000e-02, 4.466408891e-04)
          (1.078125000e-02, 1.083516120e-04)
      };
      \label{pgfplots:p1}
      \addplot [p2
      ] coordinates {
          (1.725000000e-01, 2.002334079e-03)
          (8.625000000e-02, 1.940661674e-04)
          (4.312500000e-02, 2.199212038e-05)
          (2.156250000e-02, 2.626249512e-06)
      };
      \label{pgfplots:p2}
      \addplot [p3
      ] coordinates {
         (1.725000000e-01, 9.587387784e-05)
         (8.625000000e-02, 4.903961794e-06)
         (4.312500000e-02, 2.661234226e-07)
      };
      \label{pgfplots:p3}
      \addplot [p4
      ] coordinates {
         (1.725000000e-01, 4.005809410e-06)
         (8.625000000e-02, 3.462848025e-07)
         (4.312500000e-02, 4.396345029e-08)
      };
      \label{pgfplots:p4}
      \end{axis}
      \end{tikzpicture}
      \caption{CG}
  \end{subfigure}
  \begin{subfigure}[b]{0.49\textwidth}
    \begin{tikzpicture}
      \begin{axis}[
      x=1.5cm,
      y=0.2cm,
      axis x line=bottom,
      axis y line=left,
      xmode=log,
      ymode=log,
      xlabel=h/L,ylabel= relative $L^2$ error,
      grid=major,
      legend pos=outer north east,
      ]
      \addplot [domain=1.0e-2:2.0e-1, samples=100, dashed
      ]{x^2};
      \addplot [domain=1.0e-2:2.0e-1, samples=100, dashed
      ]{0.2*x^3};
      \addplot [domain=2.0e-2:2.0e-1, samples=100, dashed
      ]{0.1*x^4};
      \addplot [domain=4.0e-2:2.0e-1, samples=100, dashed
      ]{0.1*x^5};
      \addplot [p1
      ] coordinates {
          (1.725000000e-01, 3.051984930e-02)
          (8.625000000e-02, 7.565773624e-03)
          (4.312500000e-02, 1.860592832e-03)
          (2.156250000e-02, 4.413588680e-04)
          (1.078125000e-02, 1.076480354e-04)
      };
      \addplot [p2
      ] coordinates {
          (1.725000000e-01, 1.691572031e-03)
          (8.625000000e-02, 1.472140035e-04)
          (4.312500000e-02, 1.583827188e-05)
          (2.156250000e-02, 1.718262754e-06)
      };
      \addplot [p3
      ] coordinates {
         (1.725000000e-01, 8.897625189e-05)
         (8.625000000e-02, 4.604262655e-06)
         (4.312500000e-02, 2.566428508e-07)
      };
      \addplot [p4
      ] coordinates {
         (1.725000000e-01, 4.772861265e-06)
         (8.625000000e-02, 4.024618410e-07)
         (4.312500000e-02, 5.539324908e-08)
      };
      \end{axis}
      \end{tikzpicture}
      \caption{DG}
  \end{subfigure}
  \caption{Convergence of 3D sphere benchmark, $p=1$(\ref{pgfplots:p1}),$2$(\ref{pgfplots:p2}),$3$(\ref{pgfplots:p3}),$4$(\ref{pgfplots:p4})}\label{fig:convergence_smooth}
\end{figure}
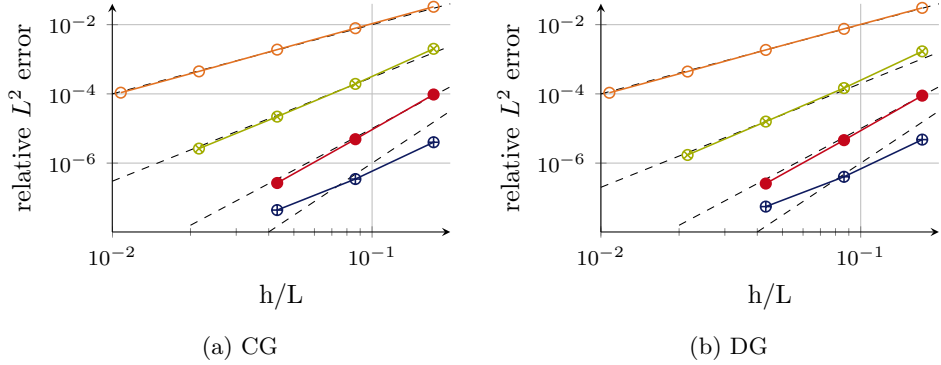
The optimal order of convergence can be observed in Figure~\ref{fig:convergence_smooth} for $p=1,2,3$ and CG and DG, whereas the slopes are not optimal for $p=4$, albeit with overall relatively low errors. The authors believe that with another choice of the volume ghost penalty parameter and well-conditioned bases for the penalty term, optimal convergence orders should also be possible for $p\geq4$.

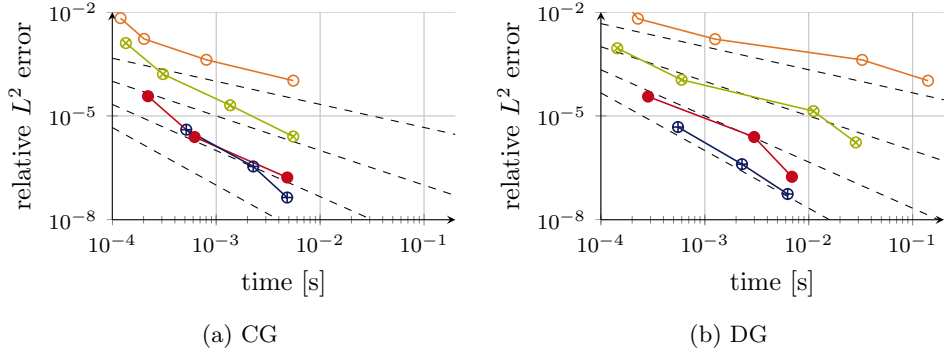
\begin{figure}
  \begin{subfigure}[b]{0.49\textwidth}
    \begin{tikzpicture}
      \begin{axis}[
      x=0.6cm,
      y=0.2cm,
      axis x line=bottom,
      axis y line=left,
      xmode=log,
      ymode=log,
      xlabel={time [s]},ylabel= relative $L^2$ error,
      grid=major,
      legend pos=outer north east,
      xmin=1e-4,
      xmax=2e-1,
      ymin=1e-8,
      ymax=1e-2,
      ]
      \addplot [domain=1.0e-4:2.0e-1, samples=100, dashed
      ]{1e-6*x^(-2/3)};
      \addplot [domain=1.0e-4:2.0e-1, samples=100, dashed
      ]{1e-8*x^(-3/3)};
      \addplot [domain=1.0e-4:2.0e-1, samples=100, dashed
      ]{1e-10*x^(-4/3)};
      \addplot [domain=1.0e-4:2.0e-1, samples=100, dashed
      ]{1e-12*x^(-5/3)};
      \addplot [p1
      ] coordinates {
          (0.0114/95, 6.695201563e-03)
          (0.02164/107, 1.677712313e-03)
          (0.08846/110, 4.186914399e-04)
          (0.5193/94, 1.045761119e-04)
      };
      \addplot [p2
      ] coordinates {
          (0.0008105/6, 1.271450330e-03)
          (0.002152/7, 1.615071096e-04)
          (0.008172/6, 2.028881766e-05)
          (0.03305/6, 2.540143154e-06)
      };
      \addplot [p3
      ] coordinates {
          (0.001324/6, 3.707483198e-05)
          (0.004962/8, 2.459981781e-06)
          (0.03862/8, 1.663371417e-07)
      };
      \addplot [p4
      ] coordinates {
          (0.004122/8, 4.005809410e-06)
          (0.06363/28, 3.462848025e-07)
          (0.09156/19, 4.396345029e-08)
      };
      \end{axis}
      \end{tikzpicture}
      \caption{CG}
  \end{subfigure}
  \begin{subfigure}[b]{0.49\textwidth}
    \begin{tikzpicture}
      \begin{axis}[
      x=0.6cm,
      y=0.2cm,
      axis x line=bottom,
      axis y line=left,
      xmode=log,
      ymode=log,
      xlabel={time [s]},ylabel= relative $L^2$ error,
      grid=major,
      legend pos=outer north east,
      xmin=1e-4,
      xmax=2e-1,
      ymin=1e-8,
      ymax=1e-2,
      ]
      \addplot [domain=1.0e-4:2.0e-1, samples=100, dashed
      ]{1e-5*x^(-2/3)};
      \addplot [domain=1.0e-4:2.0e-1, samples=100, dashed
      ]{1e-7*x^(-3/3)};
      \addplot [domain=1.0e-4:2.0e-1, samples=100, dashed
      ]{1e-9*x^(-4/3)};
      \addplot [domain=1.0e-4:2.0e-1, samples=100, dashed
      ]{1e-11*x^(-5/3)};
      \addplot [p1
      ] coordinates {
          (0.0006478/4, 2.401101056e-02)
          (0.0009081/4, 6.449337758e-03)
          (0.005024/4, 1.653418046e-03)
          (0.1307/4, 4.155972667e-04)
          (0.4215/3, 1.041812103e-04)
      };
      \addplot [p2
      ] coordinates {
        (0.0007192/5, 9.027891975e-04)
        (0.00359/6, 1.105596170e-04)
        (0.05555/5, 1.375118805e-05)
        (0.142/5, 1.714589151e-06)
      };
      \addplot [p3
      ] coordinates {
        (0.001423/5, 3.628257839e-05)
        (0.06266/21, 2.460348320e-06)
        (0.1447/21, 1.742677018e-07)
      };
      \addplot [p4
      ] coordinates {
        (0.0116/21, 4.772861278e-06)
        (0.5104/224, 4.024618269e-07)
        (0.6583/105, 5.539312170e-08)
      };
      \end{axis}
      \end{tikzpicture}
      \caption{DG}
  \end{subfigure}
  \caption{Error vs.~run time of one operator evaluation (matrix-free) of 3D sphere benchmark, $p=1$(\ref{pgfplots:p1}),$2$(\ref{pgfplots:p2}),$3$(\ref{pgfplots:p3}),$4$(\ref{pgfplots:p4})}\label{fig:error_vs_time}
\end{figure}

With the assumption that an effective preconditioner can be implemented for unfitted FEM methods with similar performance as the operator evaluation like it is the case for fitted FEM methods~\cite{Kronbichler2018} (a topic of future work), we look at the error vs.~the cost of a single operator evaluation (see Figure~\ref{fig:error_vs_time}). The results demonstrate that, under the assumption of nearly constant iteration counts as a function of the polynomial degree, a clear benefit of high-order unfitted methods with respect to runtime is possible. For a certain target accuracy, magnitudes of runtime can be saved using high-order methods.

\section{Performance}\label{section:performance}

In this section, we characterize the performance properties of the proposed algorithms by a few sample benchmarks.

\subsection{Matrix-free vs.~matrix-based method}

To compare the operator evaluation described in Sec.~\ref{section:matrixfree}, we use a fair measure for speed called \textbf{throughput} defined as
\begin{equation}
  \mathrm{throughput}=\frac{\mathrm{problem\;size\ [degrees\;of\;freedom,\;DoF]}}{\mathrm{computational\;time}\ [s]}.
\end{equation}

In the following, we characterize operator evaluation with respect to memory transfer and arithmetic operations for a benchmark on a 3D Cartesian grid resembling the background grid typically used in unfitted FEM. We compare a sparse matrix with matrix-free evaluation with the structured and unstructured quadrature algorithm (with the same structured quadrature) for a continuous (CG) and a discontinuous (DG) finite element discretization. The model equation for the benchmark is the Poisson equation~\eqref{eq:poisson} in weak form with the cell term for CG and cell and SIPG face terms for DG.

The various approaches of operator evaluation have different computational characteristics. Matrix-based methods are limited by memory bandwidth, structured matrix-free evaluation has a higher but almost constant arithmetic intensity and for unstructured matrix-free evaluation, the arithmetic intensity changes significantly from memory bound for $p=1$ to compute bound for $p\geq 3$ (discussed later together with Figure~\ref{fig:roofline}). This is challenging from an algorithm optimization perspective since overhead in floating point operations, memory access and integer overhead affect performance in the different regimes. To understand the hardware bottlenecks and identify possibilities for algorithm improvement, we estimate the arithmetic and memory complexity.

\subsubsection*{Arithmetic complexity of matrix-free evaluation}
First, we analyze the arithmetic complexity of the various approaches depending on the number of quadrature points per cell $n_\mathrm{q}^{\mathrm{dim}}$ ($n_\mathrm{q}$ denoting the one-dimensional number of quadrature points used to construct tensor-product quadratures, $n_\mathrm{q}^{\mathrm{dim}}=n_\mathrm{q}^d$) and number of DoFs in one dimension $k = p + 1$. For simplicity, we show the numbers for 3D only; the derivation for 2D is analogous.

The number of arithmetic operations (fused multiply-add operations) for evaluating the values and gradients at all quadrature points with structured tensor-product quadrature is
\begin{equation}
  \underbrace{k^3n_\mathrm{q} + k^2n_\mathrm{q}^2+n_\mathrm{q}^3k}_{\mathrm{values}}+\underbrace{3n_\mathrm{q}^4}_{\mathrm{gradients}}  \overset{n_\mathrm{q}=k}{=} 6k^4,
\end{equation}
using the collocation derivative technique presented e.g.~in~\cite{Kronbichler2019}, see also \cite[Sec.~5]{Fischer2020}. First, the DoF values are interpolated to quadrature points (collocation basis), then 1D derivatives are computed for each coordinate direction in the collocation basis.

For faces with structured tensor-product quadrature, evaluating the values and gradients involves
\begin{equation}
  \underbrace{k^2n_\mathrm{q}+n_\mathrm{q}^2k+2n_\mathrm{q}^3}_{\mathrm{values} + \mathrm{tangential\;derivatives}} + \underbrace{k^3+k^2n_\mathrm{q}+n_\mathrm{q}^2k}_{\mathrm{normal\;derivatives}} \overset{n_\mathrm{q}=k}{=} 7k^3
\end{equation}
arithmetic operations, using the collocation derivative for the tangential derivative and interpolation into the face DoFs of the normal derivative with an in-face interpolation at the face quadrature points.

For unstructured cell quadrature, evaluating all values and gradients at all quadrature involves
\begin{equation}
  \underbrace{n_\mathrm{q}^{3\mathrm{D}}(k^3 + k^2 + k)}_{\mathrm{values}} + \underbrace{n_\mathrm{q}^{3\mathrm{D}}(k^3 + 2k^2 + 3k)}_{\mathrm{gradients}} \overset{n_\mathrm{q}^{3\mathrm{D}}=k^3}{=} 2k^6+3k^5+4k^4
\end{equation}
arithmetic operations. The numbers aim to provide optimal arithmetic work by
exploiting re-occuring terms: the $d$-interpolation in $x$-direction can be
used for the partial derivatives in $y$- and $z$-directions and the
$(d-1)$-interpolation in $y$-direction can be used for partial derivative
$\partial_z$, compared to the naive interpolation of gradients for a
quadrature point ($3k^3 + 3k^2 + 3k$), cf.~\cite[Sec.~2.4]{Kronbichler2012}.

For an unstructured face for values and gradients, we have the cost
\begin{equation}
  \underbrace{n_\mathrm{q}^{2\mathrm{D}}(2k^2 + 3k)}_{\mathrm{values}+\mathrm{tangential\;derivatives}} + \underbrace{k^3 + n_\mathrm{q}^{2\mathrm{D}}(k^2 + k)}_{\mathrm{normal\;derivatives}} \overset{n_\mathrm{q}^{2\mathrm{D}}=k^2}{=} 3k^4 + 5k^3.
\end{equation}

On cut cells, depending on the shape of the interface, oftentimes refinements of cut cells are necessary for the quadrature generation~\cite{Saye2015}. This leads to an increase in the number of quadrature points in those particular cells and their intersected faces. With uniform refinement, the number of quadrature points per cell $n_\mathrm{q}^{\mathrm{dim}}$ scales with $2^{d n_\mathrm{ref}} n_\mathrm{q}^d$ where $n_\mathrm{ref}$ is the number of uniform refinements.

\subsubsection*{Memory transfer of matrix-free evaluation}
The following estimates of memory transfer are based on the actual implementation in the finite element library {deal.II}, in particular with respect to the number of pointers that are accessed per cell; for precise numbers, we refer the interested reader directly to the source code~\cite{Arndt2023}.

We estimate the memory transfer per DoF from main memory of the structured quadrature algorithm and CG with
\begin{equation}
  \frac{1}{p^d}\underbrace{\frac{14}{n_\mathrm{lanes}}}_\mathrm{Byte\;per\;cell}+\underbrace{\frac{k^d}{p^d}\,4 + 24}_\mathrm{Byte\;read/write\;DoF},
\end{equation}
where 14 bytes refer to pointers and indices to the compressed geometry data. The vectorization over cells reduces the data transfer per DoF by simultaneously setting indices and pointers for several cells. We assume that the 1D shape functions at quadrature points, as well as the weights and the geometry data, are accessible from caches. We normalize the bytes per cell by the unique DoFs per cell ($p^d$ for CG, $k^d$ for DG). Finally, 24 bytes of data transfer are necessary to read and write DoFs from and into the global vector plus 4 bytes per DoF index, which are loaded element-by-element (meaning the DoF index is duplicated for neighboring cells that share DoFs).

For DG we additionally need $2d$ face terms
\begin{equation}
  \frac{1}{k^d}\underbrace{\left(\frac{14}{n_\mathrm{lanes}}+2d\left(4+\frac{18}{n_\mathrm{lanes}}\right)\right)}_\mathrm{Byte\;per\;cell\;and\;faces}
  +\underbrace{24}_\mathrm{Byte\;read/write\;DoF},
\end{equation}
where 4 bytes of indices are necessary for each face and 18 bytes per face batch for pointers to the compressed geometry data. Again, the vectorization across faces reduces the memory transfer necessary for indices and pointers. For DG only one DoF index per cell or face is loaded while the other DoF indices can be reconstructed from the cell/face DoF index, see~\cite{Kronbichler2019}. Because even for linear elements in 3D this involves only $\frac{4}{8}$ bytes, we omit the contribution of the DoF indices in our memory transfer estimate for DG.

For unstructured quadrature, the algorithm is used together with the structured infrastructure, replacing the evaluation and integration. We switch from vectorization across cells to vectorization across quadrature points for the unstructured quadrature algorithm, since the number of quadrature points might differ for each cell. As a consequence, pointers and indices have to be set for each cell/face separately, losing the desirable property of setting them once per batch. Furthermore, our design of vectorization across quadrature points means that always the full data of a batch of points has to be loaded from memory (to avoid masked loads for the general case), even if the number of quadrature points in a batch does not fill all lanes $n_\mathrm{lanes}$. To consider this effect in our estimates we introduce $\tilde{n}_\mathrm{q} = \lceil n_\mathrm{q}/n_\mathrm{lanes} \rceil n_\mathrm{lanes}$ which is the number of quadrature points rounded up to a multiple of $n_\mathrm{lanes}$. For CG we can then estimate
\begin{equation}
  \frac{1}{p^d}\underbrace{\left(25 + \frac{14}{n_\mathrm{lanes}}\right)}_\mathrm{Byte\;per\;cell}+\frac{\tilde{n}_\mathrm{q}^d}{p^d}\underbrace{32}_\mathrm{Byte\;per\;quadrature\;point}+\underbrace{\frac{k^d}{p^d}\,4 + 24}_\mathrm{Byte\;read/write\;DoF}.
\end{equation}
Additionally to the part from the structured quadrature algorithm, we have 25 bytes per cell to set pointers to geometry data and 32 bytes per quadrature point, because its position (24 bytes) and the Jacobian determinant multiplied by the quadrature weight (\texttt{JxW}, 8 bytes) have to be loaded from main memory in contrast to the structured quadrature algorithm.

For DG and unstructured quadrature, we additionally get contributions from $2d$ faces
\begin{equation}
  \begin{split}
    &\frac{1}{k^d}\underbrace{\left[\left(25 + \frac{14}{n_\mathrm{lanes}}\right)+2d\left(29 + \frac{18}{n_\mathrm{lanes}}\right)\right]}_\mathrm{Byte\;per\;cell\;and\;faces}\\
    &+\frac{\tilde{n}_\mathrm{q}^d}{k^d}\underbrace{32}_\mathrm{Byte\;per\;cell\;qp}
    +2d\frac{\tilde{n}_\mathrm{q}^{d-1}}{k^d}\underbrace{24}_\mathrm{Byte\;per\;face\;qp}
    +\underbrace{24}_\mathrm{Byte\;read/write\;DoF}.
  \end{split}
\end{equation}
This implies 29 bytes of pointers to geometry data and 24 bytes for the $d-1$ position of the quadrature points in the face (16 bytes) plus \texttt{JxW} (8 bytes).

\subsubsection*{Memory transfer and arithmetic intensity of sparse matrix-vector multiplication}
Memory transfer and arithmetic operations per unknown of a sparse matrix-vector multiplication are of order $\mathcal{O}(p^d)$ for CG and $\mathcal{O}(k^d)$ for DG. Hence, the arithmetic intensity remains constant for sparse matrices (see Figure~\ref{fig:characteristic}). One fused multiply-add (FMA) is performed per entry, thus the intensity is approximately $\frac{2\mathrm{FLOPs}}{(8 + 4) \mathrm{Byte}} = \frac{1}{6} \frac{\mathrm{FLOPs}}{\mathrm{Byte}}$. This means the algorithm is severely limited by memory bandwidth on current computer hardware, where the machine balance indicating the hardware's capability is often around 5 -- 20\,$\frac{\mathrm{FLOPs}}{\mathrm{Byte}}$.

In detail, the memory transfer of a sparse matrix for CG scales with
$12(p+2)^d\;\mathrm{Byte}$
per DoF, assuming the column indices in the CRS format are stored as 32-bit integer numbers and that the row-start indices are negligible in size. For DG, the transfer per DoF is higher with approximately
$12(2d+1)(p+1)^d\;\mathrm{Byte}$
because of the additional face coupling matrices on the off-diagonals. We recall the approximate intensity for sparse matrices of $\frac{1}{6} \frac{\mathrm{FLOPs}}{\mathrm{Byte}}$ and model the FLOPs per DoF with
$2(p+2)^d\;\mathrm{FLOPs}$
for CG and
$2(2d+1)(p+1)^d\;\mathrm{FLOPs}$
for DG. For both CG and DG, the memory complexity of sparse matrices increases significantly with $d$ and $p$, especially for $d=3$ and $p>1$. As the application of sparse matrices is memory-bound, this leads to low throughput (discussed later together with Figure~\ref{fig:throughput}).

\subsubsection*{Verification of estimates for memory and arithmetic complexity by measurements}
The resulting characteristics from the estimates are shown in supplement~SM3 for $p=1,2,3,4$. To verify the proposed formulas for memory and arithmetic intensity, we measured the bytes and FLOPs per DoF with the performance monitoring tool \texttt{LIKWID}~\cite{Treibig2010}. The estimates fit very well to the measurements, see Figure~\ref{fig:characteristic}.

\pgfplotsset{sparsematrix/.style={tumorange,mark=otimes, semithick}}
\pgfplotsset{structured/.style={diag_violet,mark=o, semithick}}
\pgfplotsset{unstructured/.style={tumblue,mark=*, semithick}}

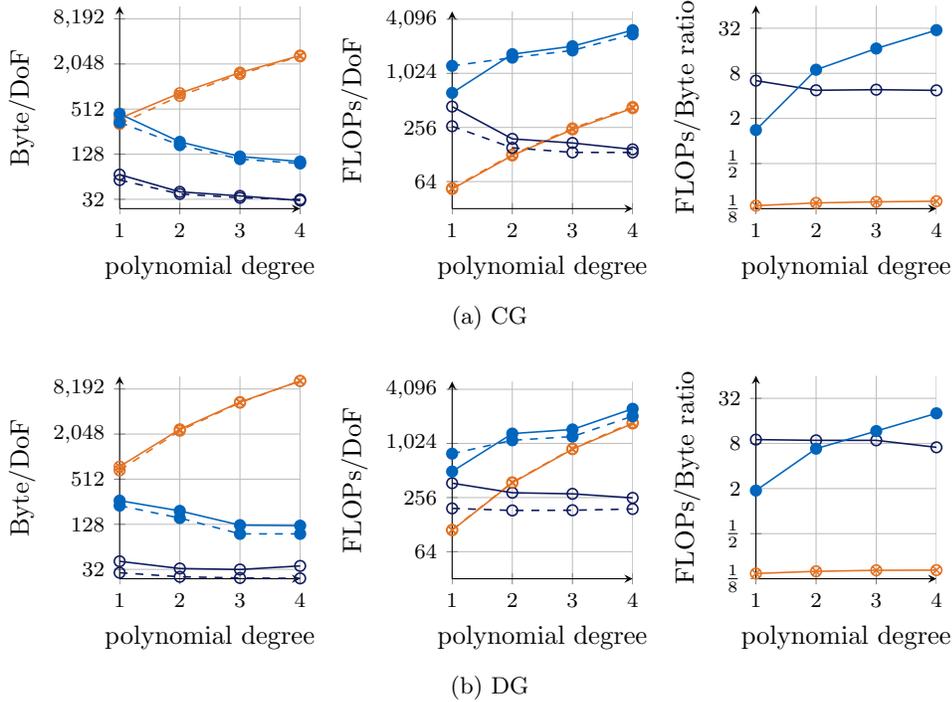
\begin{figure}
  \centering
  \begin{subfigure}[b]{\textwidth}
    \begin{tikzpicture}
      \begin{axis}[
      x=0.8cm,
      y=0.3cm,
      axis x line=bottom,
      axis y line=left,
      ymode=log,
      log basis y={2},
      ymax=12000, ymin=24,
      ytick={32, 128, 512, 2048, 8192},
      yticklabels={32, 128, 512, {2,048}, {8,192}},
      tick label style={font=\footnotesize},
      xlabel=polynomial degree,ylabel=Byte/DoF,
      log number format basis/.code 2 args={$\pgfmathparse{#1^(#2)}\pgfmathprintnumber{\pgfmathresult}$},
      grid=major,
      legend pos=outer north east,
      ]
      \addplot [sparsematrix
      ] coordinates {
          (1, 384.8) (2, 835.2) (3, 1570.5) (4, 2655.5)
      };
      \addplot [structured
      ] coordinates {
          (1, 68.3) (2, 40.2) (3, 35.5) (4, 30.8)
      };
      \addplot [unstructured
      ] coordinates {
          (1, 442) (2, 187) (3, 119) (4, 102)
      };
      \pgfplotstableread{scripts/sparse_matrix_cg.csv}\loadedtable
      \pgfplotstabletranspose\tloadedtable{\loadedtable}
      \addplot[sparsematrix, dashed, mark options={solid}] table[x expr=\coordindex+1, y index=1]{\tloadedtable};
      \pgfplotstableread{scripts/structured_cg.csv}\loadedtable
      \pgfplotstabletranspose\tloadedtable{\loadedtable}
      \addplot[structured, dashed, mark options={solid}] table[x expr=\coordindex+1, y index=1]{\tloadedtable};
      \pgfplotstableread{scripts/unstructured_cg.csv}\loadedtable
      \pgfplotstabletranspose\tloadedtable{\loadedtable}
      \addplot[unstructured, dashed, mark options={solid}] table[x expr=\coordindex+1, y index=1]{\tloadedtable};
      \end{axis}
      \end{tikzpicture}
      \begin{tikzpicture}
        \begin{axis}[
        x=0.8cm,
        y=0.36cm,
        axis x line=bottom,
        axis y line=left,
        ymode=log,
        log basis y={2},
        ymax=5000, ymin=32,
        ytick={64, 256, 1024, 4096},
        yticklabels={64, 256, {1,024}, {4,096}},
        tick label style={font=\footnotesize},
        xlabel=polynomial degree,ylabel=FLOPs/DoF,
        log number format basis/.code 2 args={$\pgfmathparse{#1^(#2)}\pgfmathprintnumber{\pgfmathresult}$},
        grid=major,
        legend pos=outer north east,
        ]
        \addplot [sparsematrix
        ] coordinates {
            (1, 52.9) (2, 124.8) (3, 243.0) (4, 419.2)
        };
        \addplot [structured
        ] coordinates {
            (1, 437) (2, 191) (3, 172) (4, 146)
        };
        \addplot [unstructured
        ] coordinates {
            (1, 620) (2, 1673) (3, 2049) (4, 3091)
        };
        \pgfplotstableread{scripts/sparse_matrix_cg.csv}\loadedtable
        \pgfplotstabletranspose\tloadedtable{\loadedtable}
        \addplot[sparsematrix, dashed, mark options={solid}] table[x expr=\coordindex+1, y index=2]{\tloadedtable};
        \pgfplotstableread{scripts/structured_cg.csv}\loadedtable
        \pgfplotstabletranspose\tloadedtable{\loadedtable}
        \addplot[structured, dashed, mark options={solid}] table[x expr=\coordindex+1, y index=4]{\tloadedtable};
        \pgfplotstableread{scripts/unstructured_cg.csv}\loadedtable
        \pgfplotstabletranspose\tloadedtable{\loadedtable}
        \addplot[unstructured, dashed, mark options={solid}] table[x expr=\coordindex+1, y index=5]{\tloadedtable};
        \end{axis}
        \end{tikzpicture}
        \begin{tikzpicture}
          \begin{axis}[
          x=0.8cm,
          y=0.3cm,
          axis x line=bottom,
          axis y line=left,
          ymode=log,
          log basis y={2},
          ymax=63.9, ymin=0.125,
          ytick={0.125, 0.5, 2, 8, 32},
          yticklabels={$\frac 18$, $\frac 12$, 2, 8, 32},
          tick label style={font=\footnotesize},
          xlabel=polynomial degree,ylabel=FLOPs/Byte ratio,
          log ticks with fixed point,
          grid=major,
          legend pos=outer north east,
          ]
          \addplot [sparsematrix
          ] coordinates {
              (1, 0.137) (2, 0.149) (3, 0.154) (4, 0.157)
          };
          \label{pgfplots:sparsematrix}
          \addplot [structured
          ] coordinates {
              (1, 6.39) (2, 4.75) (3, 4.84) (4, 4.74)
          };
          \label{pgfplots:structured}
          \addplot [unstructured
          ] coordinates {
              (1, 1.40) (2, 8.93) (3, 17.16) (4, 30.189)
          };
          \label{pgfplots:unstructured}
          \end{axis}
          \end{tikzpicture}
          \caption{CG}
  \end{subfigure}
  \begin{subfigure}[b]{\textwidth}
    \begin{tikzpicture}
      \begin{axis}[
      x=0.8cm,
      y=0.3cm,
      axis x line=bottom,
      axis y line=left,
      ymode=log,
      log basis y={2},
      ytick={32, 128, 512, 2048, 8192},
      yticklabels={32, 128, 512, {2,048}, {8,192}},
      tick label style={font=\footnotesize},
      ymax=12000, ymin=24,
      xlabel=polynomial degree,ylabel=Byte/DoF,
      log number format basis/.code 2 args={$\pgfmathparse{#1^(#2)}\pgfmathprintnumber{\pgfmathresult}$},
      grid=major,
      legend pos=outer north east,
      ]
      \addplot [sparsematrix
      ] coordinates {
          (1, 754) (2, 2369) (3, 5466) (4, 10484)
      };
      \addplot [structured
      ] coordinates {
          (1, 41) (2, 32.9) (3, 31.9) (4, 35.6)
      };
      \addplot [tumblue,mark=*, semithick
      ] coordinates {
          (1, 264.8) (2, 193.3) (3, 124.6) (4, 122.8)
      };
      \pgfplotstableread{scripts/sparse_matrix_dg.csv}\loadedtable
      \pgfplotstabletranspose\tloadedtable{\loadedtable}
      \addplot[sparsematrix, dashed, mark options={solid}] table[x expr=\coordindex+1, y index=1]{\tloadedtable};
      \pgfplotstableread{scripts/structured_dg.csv}\loadedtable
      \pgfplotstabletranspose\tloadedtable{\loadedtable}
      \addplot[structured, dashed, mark options={solid}] table[x expr=\coordindex+1, y index=1]{\tloadedtable};
      \pgfplotstableread{scripts/unstructured_dg.csv}\loadedtable
      \pgfplotstabletranspose\tloadedtable{\loadedtable}
      \addplot[unstructured, dashed, mark options={solid}] table[x expr=\coordindex+1, y index=1]{\tloadedtable};
      \end{axis}
      \end{tikzpicture}
      \begin{tikzpicture}
        \begin{axis}[
        x=0.8cm,
        y=0.36cm,
        axis x line=bottom,
        axis y line=left,
        ymode=log,
        log basis y={2},
        ytick={64, 256, 1024, 4096},
        yticklabels={64, 256, {1,024}, {4,096}},
        tick label style={font=\footnotesize},
        ymax=5000, ymin=32,
        xlabel=polynomial degree,ylabel=FLOPs/DoF,
        log number format basis/.code 2 args={$\pgfmathparse{#1^(#2)}\pgfmathprintnumber{\pgfmathresult}$},
        grid=major,
        legend pos=outer north east,
        ]
        \addplot [sparsematrix
        ] coordinates {
           (1, 111) (2, 373) (3, 884) (4, 1706)
        };
        \addplot [structured
        ] coordinates {
            (1, 370.0) (2, 289.7) (3, 281.6) (4, 253.9)
        };
        \addplot [tumblue,mark=*, semithick
        ] coordinates {
            (1, 499) (2, 1317.9) (3, 1464.5) (4, 2487)
        };
        \pgfplotstableread{scripts/sparse_matrix_dg.csv}\loadedtable
        \pgfplotstabletranspose\tloadedtable{\loadedtable}
        \addplot[sparsematrix, dashed, mark options={solid}] table[x expr=\coordindex+1, y index=2]{\tloadedtable};
        \pgfplotstableread{scripts/structured_dg.csv}\loadedtable
        \pgfplotstabletranspose\tloadedtable{\loadedtable}
        \addplot[structured, dashed, mark options={solid}] table[x expr=\coordindex+1, y index=4]{\tloadedtable};
        \pgfplotstableread{scripts/unstructured_dg.csv}\loadedtable
        \pgfplotstabletranspose\tloadedtable{\loadedtable}
        \addplot[unstructured, dashed, mark options={solid}] table[x expr=\coordindex+1, y index=5]{\tloadedtable};
        \end{axis}
        \end{tikzpicture}
        \begin{tikzpicture}
          \begin{axis}[
          x=0.8cm,
          y=0.3cm,
          axis x line=bottom,
          axis y line=left,
          ymode=log,
          log basis y={2},
          ytick={0.125, 0.5, 2, 8, 32},
          yticklabels={$\frac 18$, $\frac 12$, 2, 8, 32},
          tick label style={font=\footnotesize},
          ymax=63.9, ymin=0.125,
          xlabel=polynomial degree,ylabel=FLOPs/Byte ratio,
          log ticks with fixed point,
          grid=major,
          legend pos=outer north east,
          ]
          \addplot [sparsematrix
          ] coordinates {
             (1, 0.147) (2, 0.157) (3, 0.162) (4, 0.163)
          };
          \addplot [structured
          ] coordinates {
              (1, 9.0) (2, 8.8) (3, 8.8) (4, 7.1)
          };
          \addplot [tumblue,mark=*, semithick
          ] coordinates {
              (1, 1.88) (2, 6.81) (3, 11.7) (4, 20.2)
          };
          \end{axis}
          \end{tikzpicture}
          \caption{DG}
  \end{subfigure}
  \caption{Measurements (full line) vs.~estimates (dashed line): Characteristics of memory transfer and arithmetic operations for sparse matrix~\ref{pgfplots:sparsematrix}, structured quadrature algorithm~\ref{pgfplots:structured} and unstructured quadrature algorithm~\ref{pgfplots:unstructured}}\label{fig:characteristic}
\end{figure}

\subsubsection*{Roofline performance model}

\begin{figure}
  \centering
  \begin{subfigure}[b]{0.49\textwidth}
    \begin{tikzpicture}
      \begin{axis}[
      x=0.53cm,
      y=0.6cm,
      axis x line=bottom,
      axis y line=left,
      ymax=5000, ymin=48,
      ymode=log,
      log basis y={2},
      xmode=log,
      log basis x={2},
      xlabel=FLOPs/Byte ratio,ylabel=GFLOPs/s,
      log number format basis/.code 2 args={$\pgfmathparse{#1^(#2)}\pgfmathprintnumber{\pgfmathresult}$},
      xtick={0.125,0.5,2,8,32},
      xticklabels={0.125,0.5,2,8,32},
      grid=major,
      ]
      \path[name path=axis1] (axis cs:0.125,0) -- (axis cs:4.7,0);
      \path[name path=axis2] (axis cs:4.7,0) -- (axis cs:48,0);
  
      \addplot [
      name path=f1,
      tumblue,
      mark=none,
      domain=0.125:4.7,
      samples=40,
      ] {680*x} [sloped] node [pos=0.5,above,fill=white,inner sep=1pt]
      {\scriptsize memory bandwidth 680 GB/s};
      \addplot [
      name path=f2,
      tumblue,
      mark=none,
      domain=4.7:48,
      samples=40,
      ] {680*4.7}[sloped] node [pos=0.42,above,fill=white,inner sep=1pt]
      {\scriptsize arithmetic peak};
      \addplot [
      thick,
      color=blue,
      fill=tumblue, 
      fill opacity=0.3
      ]
      fill between[
        of=f1 and axis1,
      ];
      \addplot [
      thick,
      color=blue,
      fill=tumblue, 
      fill opacity=0.3
      ]
      fill between[
        of=f2 and axis2,
      ];
      \addplot [sparsematrix, only marks
      ] coordinates {
          (0.137, 79) (0.149, 86) (0.154, 89)
      };
      \addplot [structured, only marks
      ] coordinates {
          (6.39, 716) (4.75, 679) (4.84, 789) (4.74, 747)
      };
      \addplot [unstructured, only marks
      ] coordinates {
          (1.40, 574) (8.93, 1224) (17.16, 1758) (30.189, 2147)
      };
      \node[tumblue] at (axis cs:1.40, 574) [anchor=south west] {1};
      \node[tumblue] at (axis cs:8.93, 1224) [anchor=south west] {2};
      \node[tumblue] at (axis cs:17.16, 1758) [anchor=south west] {3};
      \node[tumblue] at (axis cs:30.189, 2147) [anchor=south west] {4};
      \end{axis}
      \end{tikzpicture}
      \caption{CG}
  \end{subfigure}
  \hfill
  \begin{subfigure}[b]{0.49\textwidth}
    \begin{tikzpicture}
      \begin{axis}[
      x=0.53cm,
      y=0.6cm,
      axis x line=bottom,
      axis y line=left,
      ymax=5000, ymin=48,
      ymode=log,
      log basis y={2},
      xmode=log,
      log basis x={2},
      xlabel=FLOPs/Byte ratio,ylabel=GFLOPs/s,
      log number format basis/.code 2 args={$\pgfmathparse{#1^(#2)}\pgfmathprintnumber{\pgfmathresult}$},
      xtick={0.125,0.5,2,8,32},
      xticklabels={0.125,0.5,2,8,32},
      grid=major,
      ]
      \path[name path=axis1] (axis cs:0.125,0) -- (axis cs:4.7,0);
      \path[name path=axis2] (axis cs:4.7,0) -- (axis cs:48,0);
  
      \addplot [
      name path=f1,
      tumblue,
      mark=none,
      domain=0.125:4.7,
      samples=40,
      ] {680*x} [sloped] node [pos=0.5,above,fill=white,inner sep=1pt]
      {\scriptsize memory bandwidth 680 GB/s};
      \addplot [
      name path=f2,
      tumblue,
      mark=none,
      domain=4.7:48,
      samples=40,
      ] {680*4.7}[sloped] node [pos=0.42,above,fill=white,inner sep=1pt]
      {\scriptsize arithmetic peak};
      \addplot [
      thick,
      color=blue,
      fill=tumblue, 
      fill opacity=0.3
      ]
      fill between[
        of=f1 and axis1,
      ];
      \addplot [
      thick,
      color=blue,
      fill=tumblue, 
      fill opacity=0.3
      ]
      fill between[
        of=f2 and axis2,
      ];
      \addplot [sparsematrix, only marks
      ] coordinates {
        (0.147, 83.5) (0.157, 86) (0.162, 98) (0.163, 100)
      };
      \label{pgfplots:sparsematrix_marks}
      \addplot [structured, only marks
      ] coordinates {
          (9.0, 841) (8.8, 829) (8.8, 955) (7.1, 864)
      };
      \label{pgfplots:structured_marks}
      \addplot [unstructured, only marks
      ] coordinates {
          (1.88, 455) (6.81, 1036) (11.7, 1544) (20.2, 1849)
      };
      \label{pgfplots:unstructured_marks}
      \node[tumblue] at (axis cs:1.88, 455) [anchor=south west] {1};
      \node[tumblue] at (axis cs:6.81, 1036) [anchor=south west] {2};
      \node[tumblue] at (axis cs:11.7, 1544) [anchor=south west] {3};
      \node[tumblue] at (axis cs:20.2, 1849) [anchor=south west] {4};
      \end{axis}
      \end{tikzpicture}
      \caption{DG}
  \end{subfigure}

  \caption{Roofline performance model: 3D, classification of sparse matrix (\ref{pgfplots:sparsematrix_marks}) compared to structured quadrature algorithm (\ref{pgfplots:structured_marks}) and unstructured quadrature algorithm (\ref{pgfplots:unstructured_marks}) for $p=1,2,3,4$ on AMD EPYC 9354}\label{fig:roofline}
\end{figure}
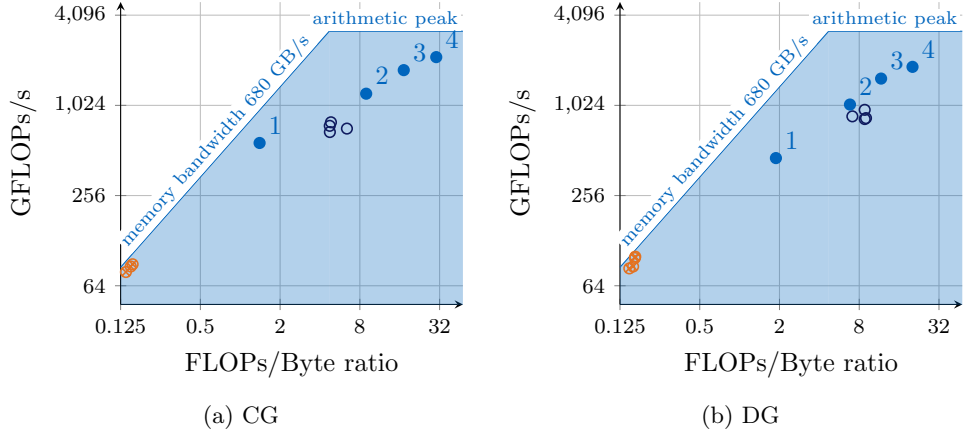

The roofline performance model~\cite{Williams2009} is a graphical representation used to assess and visualize the performance efficiency of a software implementation of an algorithm on a particular computational system by comparing achieved performance against hardware-imposed limitations.

According to Figure~\ref{fig:roofline}, both matrix-free algorithms reach a substantially higher arithmetic intensity than the sparse matrix. The intensity stays nearly constant for the structured quadrature, whereas the arithmetic intensity of the unstructured algorithm increases with the polynomial degree. The proximity to the different hardware bottlenecks, namely the main memory bandwidth and the arithmetic peak performance, shows the algorithms' high optimization level. For $p=4$, the unstructured algorithm reaches 65\% of peak arithmetic performance for CG and 55\% for DG, respectively.
This is less than the values of 85--90\% typically reached for dense matrix-matrix multiplication (\texttt{dgemm})~\cite{VanZee2015}; the main reason for lower performance is the unstructured read/write access into the global vectors and operations at quadrature points, which lead to non-floating-point instructions and affect the instruction pipelining.

The position in the roofline plot does not reveal the actual speed of operator evaluation. To evaluate the performance of the various algorithms, we next assess the throughput of the operator evaluation.

\subsubsection*{Throughput}

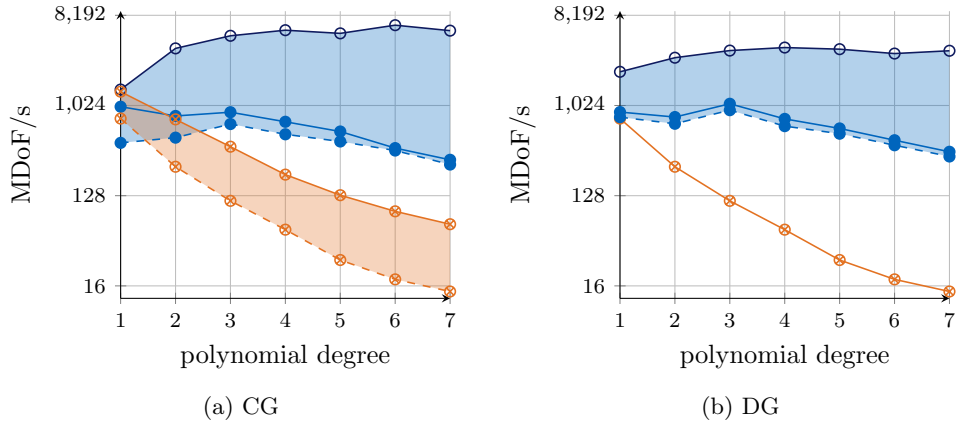
\begin{figure}
  \centering
  \begin{subfigure}[b]{0.49\textwidth}
    \begin{tikzpicture}
      \begin{axis}[
      x=0.73cm,
      y=0.4cm,
      axis x line=bottom,
      axis y line=left,
      ymode=log,
      log basis y={2},
      ymax=9000, ymin=12,
      xlabel=polynomial degree,ylabel=MDoF/s,
      log number format basis/.code 2 args={$\pgfmathparse{#1^(#2)}\pgfmathprintnumber{\pgfmathresult}$},
      xtick={1,...,7},
      grid=major,
      legend pos=outer north east,
      ]
      \addplot [name path=upper_mf, diag_violet,mark=o, semithick
      ] coordinates {
          (1, 1477.4) (2, 3815.92) (3, 5102.54) (4, 5804.79) (5, 5395.06) (6, 6515.98) (7, 5744.42)
      };
      \addplot [tumblue,mark=*, semithick
      ] coordinates {
          (1, 994.763) (2, 803.726) (3, 876.202) (4, 702.942) (5, 561.621) (6, 383.209) (7, 292.63)
      };
      \addplot [name path=lower_mf, tumblue, mark=*, dashed, semithick, mark options={solid}
      ] coordinates {
          (1, 432.219) (2, 487.474) (3, 667.124) (4, 526.496) (5, 446.192) (6, 362.28) (7, 262.203)
      };
      \addplot [name path=upper_ma, tumorange, mark=otimes, semithick
      ] coordinates {
          (1, 1406.04) (2, 743.003) (3, 395.838) (4, 207.926) (5, 129.366) (6, 89.1993) (7, 66.2796)
      };
      \addplot [name path=lower_ma, tumorange, mark=otimes, dashed, semithick, mark options={solid}
      ] coordinates {
          (1, 756.805) (2, 249.274) (3, 113.586) (4, 58.4536) (5, 29.0362) (6, 18.586) (7, 14.0511)
      };
      \label{pgfplots:sparsematrix_lower_bound}
      \addplot [
        thick,
        color=tumblue,
        fill=tumblue, 
        fill opacity=0.3
        ]
        fill between[
          of=upper_mf and lower_mf,
        ];
      \addplot [
        thick,
        color=tumorange,
        fill=tumorange, 
        fill opacity=0.3
        ]
        fill between[
          of=upper_ma and lower_ma,
        ];
        \label{pgfplots:sparsematrix_range}
      \end{axis}
      \end{tikzpicture}
      \caption{CG}
  \end{subfigure}
  \hfill
  \begin{subfigure}[b]{0.49\textwidth}
    \begin{tikzpicture}
      \begin{axis}[
      x=0.73cm,
      y=0.4cm,
      axis x line=bottom,
      axis y line=left,
      ymode=log,
      log basis y={2},
      ymax=9000, ymin=12,
      xlabel=polynomial degree,ylabel=MDoF/s,
      log number format basis/.code 2 args={$\pgfmathparse{#1^(#2)}\pgfmathprintnumber{\pgfmathresult}$},
      xtick={1,...,7},
      grid=major,
      legend pos=outer north east,
      ]
      \addplot [name path=upper_mf, diag_violet,mark=o, semithick
      ] coordinates {
          (1, 2224.15) (2, 3078.61) (3, 3629.1) (4, 3891.29) (5, 3747.89) (6, 3387.17) (7, 3612.92)
      };
      \addplot [tumblue,mark=*, semithick
      ] coordinates {
          (1, 879.436) (2, 783.792) (3, 1064.25) (4, 751.301) (5, 602.738) (6, 459.27) (7, 352.02)
      };
      \addplot [tumorange,mark=otimes, semithick
      ] coordinates {
          (1, 756.805) (2, 249.274) (3, 113.586) (4, 58.4536) (5, 29.0362) (6, 18.586) (7, 14.0511)
      };
      \addplot [name path=lower_mf, tumblue, mark=*, dashed, semithick, mark options={solid}
      ] coordinates {
          (1, 781.738) (2, 671.42) (3, 918.249) (4, 635.698) (5, 534.012) (6, 409.132) (7, 315.801)
      };
      \label{pgfplots:unstructured_with_surface}
      \addplot [
      thick,
      color=blue,
      fill=tumblue, 
      fill opacity=0.3
      ]
      fill between[
        of=upper_mf and lower_mf,
      ];
      \label{pgfplots:matrixfree_range}
      \end{axis}
      \end{tikzpicture}
      \caption{DG}
  \end{subfigure}

  \caption{Throughput: 3D, sparse matrix upper bound (\ref{pgfplots:sparsematrix}) and sparse matrix lower bound (\ref{pgfplots:sparsematrix_lower_bound}) enclosing the sparse matrix range (\ref{pgfplots:sparsematrix_range}) compared to structured quadrature algorithm (\ref{pgfplots:structured}), unstructured quadrature algorithm with standard terms (\ref{pgfplots:unstructured}) and with additional surface terms (\ref{pgfplots:unstructured_with_surface}) enclosing the matrix-free range (\ref{pgfplots:matrixfree_range}) for $p=1,\ldots,7$ on AMD EPYC 9354.}\label{fig:throughput}
\end{figure}

Figure~\ref{fig:throughput} shows the achievable throughput range of sparse matrix vs.~matrix-free operator evaluation. The range bounds depend on the cut-cell ratio. For the sparse matrix and CG the worst-case scenario (lower bound) is that all faces need ghost penalty stabilization which results in a sparsity pattern (and memory complexity) like for DG. For the matrix-free evaluation, the lower bound represents unstructured quadrature evaluation on all cells with an additional surface term (representing the Nitsche term~\eqref{eq:bilinear_cg}). We observe a higher throughput for the proposed fully matrix-free algorithm compared to the matrix-based algorithm for CG for $p\geq3$ and for DG for all degrees on the EPYC 9354 machine. As systems with higher arithmetic performance compared to memory bandwidth are more common than the chosen system on recent HPC installations, we expect that matrix-free methods outperform sparse matrices even more substantially.

From the complexity analysis, we would expect to obtain the highest performance for unstructured quadrature for linear elements. One reason for why this is not the case is that SIMD vectorization across quadrature points does not fit to the number of quadrature points of surface and face integrals for $p=1,2$. Assuming $n_\mathrm{lanes}=8$ and $n_{\mathrm{q,face/surface}}=(p+1)^2\overset{p=1}{=}4$ for linear elements, the floating point operations and memory transfer are doubled. To overcome this problem, future code could adaptively select precompiled code paths for available SIMD vectorization widths (scalar to maximal available width $n_\mathrm{lanes}$). The same reasoning applies to $p=2$. Furthermore, the linear code path uses additional specialized functions that directly inject the properties of the polynomial space, whereas all code paths with $p\geq 2$ suffer from suboptimal code generation due to loop lengths to be determined during runtime.

\begin{remark}
  The presented numbers are based on the operator evaluation only. When also taking the time to compute and assemble the global matrix into account, the proposed matrix-free algorithms outperform a sparse matrix-based solver when looking at total wall time even if the matrix-free operator evaluation is slightly slower than matrix-vector multiplication.
\end{remark}

As the matrix-free algorithms still have an overhead, indicated by a gap to the performance ceiling in the roofline model especially for low polynomial degrees, it is likely that further code optimizations could accelerate the matrix-free performance even more, outperforming sparse matrices for any combination of dimension, degree, and type of finite element discretization for cheap operators like arising from the Poisson equation. Furthermore, the performance models indicate that implementations for GPUs~\cite{Chalmers2023, Kolev2021, KronbichlerLjungkvist2019, Ljungkvist2017} would show similar performance advantages.

\subsubsection*{Unstructured quadrature terms vs.~structured quadrature stabilization terms}

Compared to standard (fitted) finite element computations, the unfitted setting involves cut elements with unstructured quadrature and stabilization terms. In terms of performance impact, the evaluation of unstructured quadrature poses a greater challenge than the stabilization terms since the latter still exposes quadrature points as a tensor product, and therefore, enables lower complexity algorithms. Especially for the DG case, flux terms on cut faces are typically the dominating factor.

\subsection{Plane benchmark}
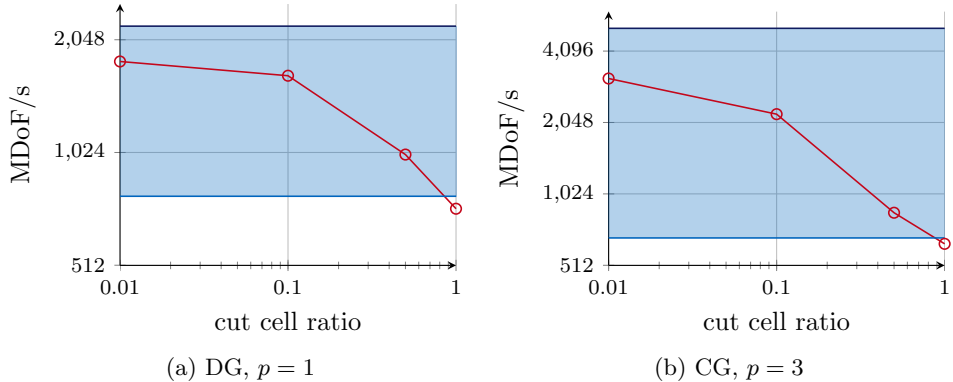
\begin{figure}
  \centering
  \begin{subfigure}[b]{0.49\textwidth}
    \begin{tikzpicture}
      \begin{axis}[
      x=0.97cm,
      y=1.5cm,
      axis x line=bottom,
      axis y line=left,
      xmode=log,
      ymode=log,
      log basis y={2},
      ymax=2548, ymin=511,
      xlabel=cut cell ratio,ylabel=MDoF/s,
      log number format basis/.code 2 args={$\pgfmathparse{#1^(#2)}\pgfmathprintnumber{\pgfmathresult}$},
      grid=major,
      legend pos=outer north east,
      ]
      \addplot [name path=upper_mf, diag_violet, semithick
      ] coordinates {
          (1, 2224.15) (0.01, 2224.15)
      };
      \addplot [tumred,mark=o, semithick
      ] coordinates {
          (1, 724) (0.5, 1010) (0.1, 1640) (0.01, 1790) (0, 1930)
      };
      \label{pgfplots:combined_ratio_cut_cells}
      \addplot [name path=lower_mf, tumblue, semithick
      ] coordinates {
          (1, 781.738) (0.01, 781.738)
      };
      \addplot [
        thick,
        color=blue,
        fill=tumblue, 
        fill opacity=0.3
        ]
        fill between[
          of=upper_mf and lower_mf,
        ];
      \end{axis}
    \end{tikzpicture}
    \caption{DG, $p=1$}
  \end{subfigure}
  \begin{subfigure}[b]{0.49\textwidth}
    \begin{tikzpicture}
      \begin{axis}[
      x=0.97cm,
      y=0.95cm,
      axis x line=bottom,
      axis y line=left,
      xmode=log,
      ymode=log,
      log basis y={2},
      ymax=6000, ymin=511,
      xlabel=cut cell ratio,ylabel=MDoF/s,
      log number format basis/.code 2 args={$\pgfmathparse{#1^(#2)}\pgfmathprintnumber{\pgfmathresult}$},
      grid=major,
      legend pos=outer north east,
      ]
      \addplot [name path=upper_mf, diag_violet, semithick
      ] coordinates {
          (1, 5102.54) (0.01, 5102.54)
      };
      \addplot [tumred,mark=o, semithick
      ] coordinates {
          (1, 630) (0.5, 852) (0.1, 2218) (0.01, 3135)
      };
      \addplot [name path=lower_mf, tumblue, semithick
      ] coordinates {
          (1, 667.124) (0.01, 667.124)
      };
      \addplot [
        thick,
        color=blue,
        fill=tumblue, 
        fill opacity=0.3
        ]
        fill between[
          of=upper_mf and lower_mf,
        ];
      \end{axis}
    \end{tikzpicture}
    \caption{CG, $p=3$}\label{fig:combined_throughput_p3}
  \end{subfigure}
  \caption{Throughput of combined matrix-free operator evaluation depending on the ratio of cut cells to total physical cells~(\ref{pgfplots:combined_ratio_cut_cells}) compared to the throughput range suggested by the generic throughput benchmark~(\ref{pgfplots:matrixfree_range} from Figure~\ref{fig:throughput})}\label{fig:combined_ratio_cut_cells}
\end{figure}

In order to quantify the performance envelope of the combined structured/unstructured matrix-free evaluation, we use a benchmark with a controlled ratio of cut elements compared to all physical elements. This is achieved by choosing a level set function as a plane with an adjustable position. In Figure~\ref{fig:combined_ratio_cut_cells}, we show the throughput of the combined matrix-free evaluation compared to the range the idealized benchmark would suggest for DG with degree $p=1$ and CG with $p=3$. When all cells are cut, we observe a minor degradation of throughput because of the additional stabilization terms, which have not been considered in the benchmark above, and because of the logic to select the correct combined code path. The dependence on the cut cell ratio is clearly visible, which shows that structured matrix-free evaluation routines with tensor products should be used whenever possible, giving a significant speedup of the overall operator evaluation. Even for a cut-cell ratio of 1, the throughput is bounded by the worst-case scenario of unstructured quadrature in all cells. For higher polynomial degrees, the speedup due to combined matrix-free evaluation is even higher, see Figure~\ref{fig:combined_throughput_p3}.

\section{Application}\label{section:application}

\subsection{Large-scale simulations}

\begin{figure}
  \centering
  \begin{subfigure}[b]{0.49\textwidth}
    \includegraphics[trim=350 80 350 150,clip,width=\textwidth]{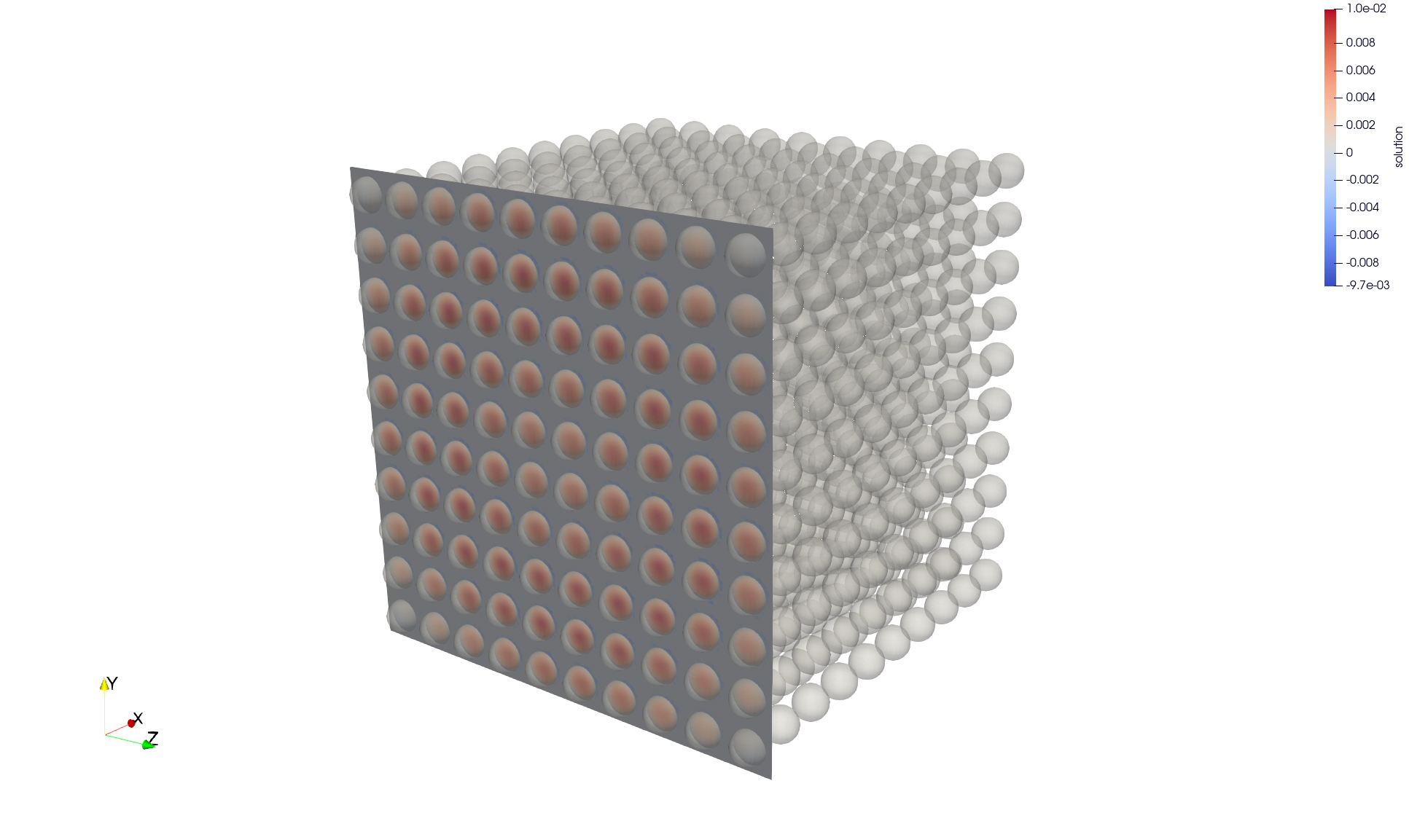}
    \caption{Regularly arranged 1000 spheres}\label{fig:spheres_regular}
  \end{subfigure}
  \begin{subfigure}[b]{0.49\textwidth}
    \includegraphics[trim=350 80 350 150,clip,width=\textwidth]{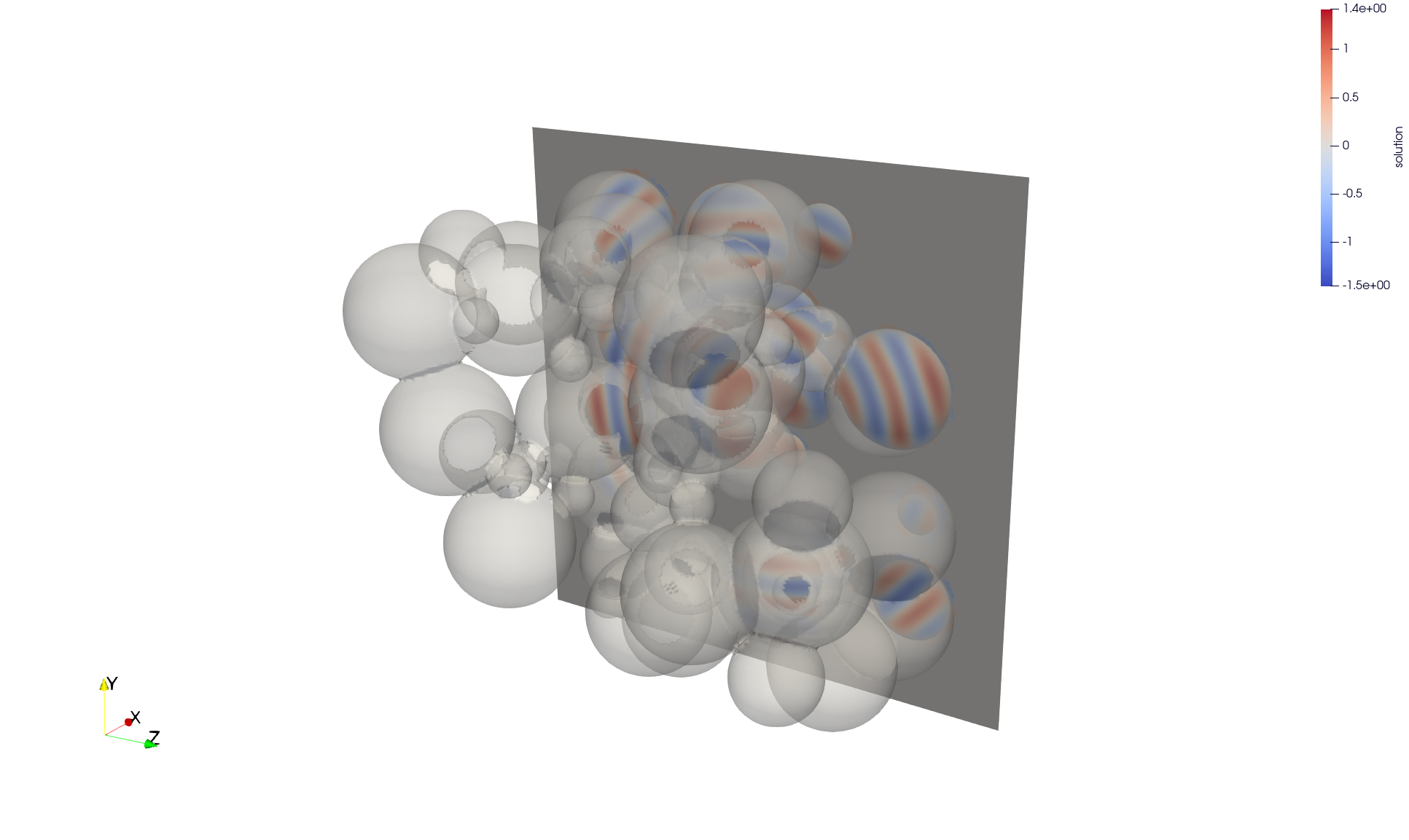}
    \caption{Random 100 spheres}\label{fig:spheres_random}
  \end{subfigure}
  \caption{Sphere benchmark cases}
\end{figure}

First, we show the operator evaluation for 1000 spheres (see Figure~\ref{fig:spheres_regular}) arranged in a regular pattern. Then, we present a computation on a domain whose surface is a random distribution of intersecting spheres (see Figure~\ref{fig:spheres_random}), which resembles the geometry of an anode of a lithium battery~\cite{Schmidt2023}. As a background grid, we use a mesh with $12^3$ cells, where we apply uniform global refinements.
The embedded geometries are described by level-set functions. The signed distance level-set function of a single sphere is defined as
\begin{equation}
  \phi(\Vector{x})=\norm{\Vector{x}-\Vector{x}_\mathrm{c}}-R,
\end{equation}
where $\Vector{x}_\mathrm{c}$ denotes the center of the sphere and $R$ the radius. The union of two spheres can be calculated by
\begin{equation}
  \phi_1 \cup \phi_2 = \min(\phi_1,\phi_2).
\end{equation}
This Boolean operation allows us to create level-set functions of arbitrarily combined intersecting spheres.

The two most important factors for throughput are the ratio of cut cells and the chosen load-balancing strategy between cut and uncut cells. For the randomly arranged spheres, we observe a number of quadrature points per cut cell of over 1000 for some cells. As the load balancing strategy solely relies on the classification of cells into cut or inside cells, the load imbalance caused by different numbers of quadrature points is not accounted for. Therefore, a more sophisticated load-balancing algorithm, considering the number of quadrature points per cell in a cost model, could further enhance the performance of matrix-free algorithms.

For the regular spheres benchmark with 15,264,256 DoF, we record a throughput of 802 MDoF/s at a cut cell ratio of 0.478 for a discontinuous Galerkin discretization with polynomial degree $p=1$.

For the random sphere benchmark with 13,139,208 DoF, we measure a throughput of 1042 MDoF/s at a cut cell ratio of 0.178 for a discontinuous Galerkin discretization with degree $p=1$. If we increase the polynomial degree to $p=3$ with the same mesh, giving a problem with 105,116,160 unknowns, we observe a slightly higher throughput of 1103 MDoF/s.

The presented examples indicate that the throughput is mainly influenced by the cut-cell ratio and less by the complexity of the cuts.

\subsection{Analysis of computational cost}

\begin{figure}
  \centering
  \begin{subfigure}[b]{0.99\textwidth}
    \begin{tikzpicture}
      \begin{axis}[
      width=0.90\textwidth,
      y=0.5cm,
      legend style={at={(0.5,1.1)}, anchor=south, legend columns=-1, font=\small},
      xbar stacked,
      ytick=data,
      symbolic y coords={matrix-free, matrix-based}]
      \addplot coordinates {
      (0,matrix-free) (2.6277,matrix-based)
      };
      \addplot coordinates {
      (0.0798,matrix-free) (0.4555,matrix-based)
      };
      \addplot coordinates {
      (30.9886,matrix-free) (27.0754,matrix-based)
      };
      \addplot coordinates {
      (1.00e+01,matrix-free) (2.08e+01,matrix-based)
      };
      \legend{assembly, operator eval., setup multigrid, apply v-cycle}
      \end{axis}
    \end{tikzpicture}
    \caption{DG, $p=3$, 15 iterations, 3,672,064 DoF}\label{fig:time_p3_dg}
  \end{subfigure}
  \begin{subfigure}[b]{0.99\textwidth}
    \begin{tikzpicture}
      \begin{axis}[
      width=0.90\textwidth,
      y=0.5cm,
      xbar stacked,
      ytick=data,
      symbolic y coords={matrix-free, matrix-based}]
      \addplot coordinates {
      (0,matrix-free) (0.8114,matrix-based)
      };
      \addplot coordinates {
      (0.0258,matrix-free) (0.0354,matrix-based)
      };
      \addplot coordinates {
      (6.1701,matrix-free) (5.7715,matrix-based)
      };
      \addplot coordinates {
      (1.91,matrix-free) (1.71,matrix-based)
      };
      \end{axis}
    \end{tikzpicture}
    \caption{DG, $p=1$, 4 iterations, 3,510,336 DoF}\label{fig:time_p1_dg}
  \end{subfigure}
  \begin{subfigure}[b]{0.99\textwidth}
    \begin{tikzpicture}
      \begin{axis}[
      width=0.90\textwidth,
      y=0.5cm,
      legend style={at={(1.3,1.0)},
      anchor=north},
      xbar stacked,
      ytick=data,
      symbolic y coords={matrix-free, matrix-based}]
      \addplot coordinates {
      (0,matrix-free) (1.8708,matrix-based)
      };
      \addplot coordinates {
      (0.0164,matrix-free) (0.0941,matrix-based)
      };
      \addplot coordinates {
      (10.0435,matrix-free) (11.7496,matrix-based)
      };
      \addplot coordinates {
      (1.97,matrix-free) (3.69,matrix-based)
      };
      \end{axis}
    \end{tikzpicture}
    \caption{CG, $p=3$, 6 iterations, 1,597,645 DoF}\label{fig:time_p3_cg}
  \end{subfigure}
  \begin{subfigure}[b]{0.99\textwidth}
    \begin{tikzpicture}
      \begin{axis}[
      width=0.90\textwidth,
      y=0.5cm,
      legend style={at={(0.5,1.1)}, anchor=south, legend columns=-1, font=\small},
      xbar stacked,
      ytick=data,
      symbolic y coords={matrix-free, matrix-based}]
      \addplot coordinates {
      (0,matrix-free) (0.8239,matrix-based)
      };
      \addplot coordinates {
      (0.4499,matrix-free) (0.2228,matrix-based)
      };
      \addplot coordinates {
      (0.8019,matrix-free) (0.2276,matrix-based)
      };
      \addplot coordinates {
      (5.23e-01,matrix-free) (5.48e-01,matrix-based)
      };
      \legend{assembly, operator eval., setup AMG, apply AMG}
      \end{axis}
    \end{tikzpicture}
    \caption{CG, $p=1$, 83 iterations, 3,507,367 DoF}\label{fig:time_p1_cg}
  \end{subfigure}
  \caption{Solver time in [$s$]. For $p=3$: $\tau_\mathrm{v}$=0.001 and for $p=1$: $\tau_\mathrm{v}=1$}\label{fig:solver_time}
\end{figure}
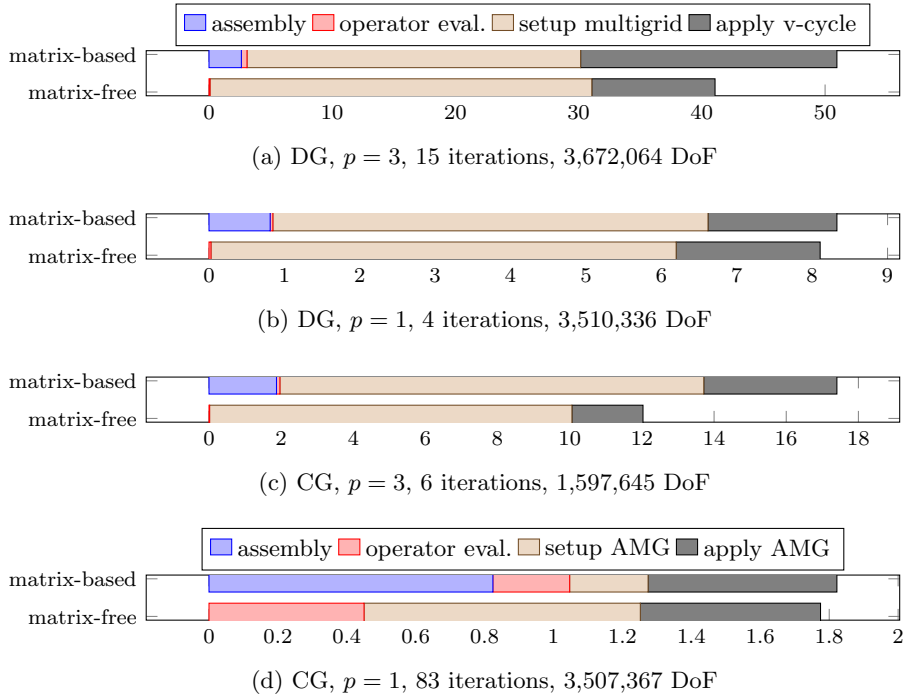

For the analysis of computational cost, we look at the single-sphere benchmark. We consider the time for the solver and neglect the setup time for e.g. repartitioning and quadrature generation because they are identical for matrix-based and matrix-free approaches. We analyze $p=1,3$ for CG and DG and compare matrix-based to matrix-free realizations in Figure~\ref{fig:solver_time}. We use the hybrid multigrid preconditioner for all examples except for CG and $p=1$, where we directly employ AMG. All simulations are conducted on a single AMD EPYC 9354 node.

We observe for all presented examples that the preconditioner consumes most of the runtime, especially for high-order methods.
While for the matrix-based methods, the assembly and the actual operator evaluation runtime are still visible in the figures, with the presented matrix-free approach the cost of operator evaluation becomes negligible. This is mainly due to the high cost of factorizing and applying the additive Schwarz smoother, which is based on matrices.
While the matrix-based preconditioner limits the speedup obtained for the chosen test case, the results are promising in the sense that much larger gains can be expected with better preconditioning.


We notice that applying the v-cycle is faster for the matrix-free path for $p=3$ (both CG and DG). This is due to the faster matrix-free operator evaluation, which is used to calculate the residual in every smoothing iteration. Additionally, the load-balancing for the Chebyshev-accelerated additive Schwarz smoother is not optimal in the matrix-based case.
Optimal load-balancing is a non-trivial topic for both matrix-free and matrix-based methods due to the conflicting requirements in different parts of the solution process. Looking at the presented graphs indicates that optimizing the load balancing for the multigrid preconditioner is currently most important.

\section{Conclusion}\label{section:conclusion}

We have discussed strategies for efficiently implementing unfitted finite element operators in a matrix-free way, utilizing current computer hardware better by avoiding sparse matrices. This leads to higher throughput, a measure for the speed of operator evaluation, by approximately an order of magnitude for $p=3$ and DG. With the proposed algorithms, significantly faster computations with high-order unfitted finite element discretizations become possible, especially in 3D.

Beyond the advancements in the present work, several paths for further improving the operator evaluation are conceivable, such as utilizing the special tensor-product structure from the quadratures generated on cut cells, implementing an adaptive SIMD strategy for increasing throughput with low polynomial degrees and reducing the integer overhead to better saturate the hardware with useful work, i.e., fetching data from memory or floating point operations.

The second important topic for future work is to develop more efficient preconditioners. While multigrid preconditioners together with additive Schwarz based smoothers have shown robust iteration counts, the factorization of local dense matrices is a serious performance bottleneck. Considerable improvements seem possible by exploring fast-diagonalization techniques combined with low-rank updates on the crucial eigenmodes, as well as basis equilibration methods to reduce ill-conditioning of basis functions evaluated in neighboring elements via the ghost penalty stabilization.


\section*{Acknowledgements}
The authors would like to thank Simon Sticko for the software contributions to the non-matching infrastructure of the deal.II finite element library, André Massing for discussions about unfitted methods and Niklas Fehn for discussions about software design.

\bibliographystyle{siamplain}
\bibliography{bibliography}
\end{document}

%% file: ex_shared.tex
\usepackage{lipsum}
\usepackage{amsfonts}
\usepackage{graphicx}
\usepackage{epstopdf}
\ifpdf
  \DeclareGraphicsExtensions{.eps,.pdf,.png,.jpg}
\else
  \DeclareGraphicsExtensions{.eps}
\fi

\newsiamremark{remark}{Remark}
\newsiamremark{hypothesis}{Hypothesis}
\crefname{hypothesis}{Hypothesis}{Hypotheses}
\newsiamthm{claim}{Claim}

\headers{Matrix-free algorithms for unfitted methods}{M. Bergbauer, P. Munch, W. A. Wall, and M. Kronbichler}

\title{High-performance matrix-free unfitted finite element operator evaluation\thanks{Submitted to the editors April 11, 2024.
\funding{We gratefully acknowledge the support from the German Research Foundation (DFG) under the project “High-Performance Cut Discontinuous Galerkin Methods for Flow Problems and Surface-Coupled Multiphysics Problems” Grant Agreement No. 456365667.}}}

\author{Maximilian Bergbauer\thanks{Technical University of Munich
  (\email{maximilian.bergbauer@tum.de}, \email{wolfgang.a.wall@tum.de}).}
\and Peter Munch\thanks{Uppsala University, University of Augsburg (\email{peter.munch@it.uu.se})}
\and Wolfgang A. Wall\footnotemark[2]
\and Martin Kronbichler \thanks{Ruhr University Bochum, University of Augsburg (\email{martin.kronbichler@rub.de}).}}

\usepackage{amsopn}

\input{DGCommands.tex}

\input{tumcolors.tex}
\usepackage{algpseudocode}
\usepackage{forest}
\usepackage{hyperref}
\usepackage{tikz}
\usetikzlibrary{arrows}
\usetikzlibrary{positioning}
\usepackage{pgfplots}
\usepgfplotslibrary{fillbetween}
\usepackage{pgfplotstable}
\pgfplotsset{compat=1.9}
\pgfplotsset{every tick label/.append style={font=\footnotesize}}
\usepackage{booktabs}
\usepackage{caption}
\usepackage[labelformat=simple]{subcaption}

\usepackage{setspace}
\reversemarginpar


%% file: DGCommands.tex
\newcommand{\Grad}[1]{\nabla{#1}}
\newcommand{\Laplace}[1]{\Delta{#1}}

\usepackage{stmaryrd}

\newcommand{\avg}[1]{\{\!\{#1\}\!\}}

\newcommand{\jump}[1]{\llbracket{#1} \rrbracket }

\newcommand{\intinteriorfaces}[2]{ \left({#1},{#2} \right)_{\Gamma_{h}^{\mathrm{int}} }}

\newcommand{\intstabilizedpatch}[2]{ \left({#1},{#2} \right)_{P}}

\newcommand{\intdom}[2]{ \left({#1},{#2} \right)_{\Omega_{h}^{\;}} }

\newcommand{\intdomDirichlet}[2]{ {\left({#1},{#2} \right)}_{\Gamma^{\rm{D}}_{h}}}

\usepackage{bm}

\newcommand{\Matrix}[1]{\bm{\mathrm{{#1}}}}
\newcommand{\Vector}[1]{\bm{\mathrm{{#1}}}}

\newcommand\norm[1]{\lVert#1\rVert}


%% file: tumcolors.tex
\definecolor{black}{RGB}{  0,   0,   0}
\definecolor{darkgray}{RGB}{ 88,  88,  90}
\definecolor{lightgray}{RGB}{217, 218, 219}
\definecolor{mediumgray}{RGB}{156, 157, 159}
\definecolor{pantone283}{RGB}{152, 198, 234}
\definecolor{pantone300}{RGB}{  0, 101, 189}
\definecolor{pantone301}{RGB}{  0,  82, 147}
\definecolor{pantone540}{RGB}{  0,  51,  89}
\definecolor{pantone542}{RGB}{100, 160, 200}
\definecolor{tumblue}{RGB}{  0, 101, 189}
\definecolor{tumdarkblue}{RGB}{  0,  82, 147}
\definecolor{tumgreen}{RGB}{162, 173,   0}
\definecolor{tumivory}{RGB}{218, 215, 203}
\definecolor{tumlightblue}{RGB}{152, 198, 234}
\definecolor{tumorange}{RGB}{227, 114,  34}
\definecolor{tumred}{RGB}{196,   7,  27}
\definecolor{white}{RGB}{255, 255, 255}

\definecolor{acc_blue}{RGB}{  0, 153, 255}
\definecolor{acc_darkred}{RGB}{202,  33,  63}
\definecolor{acc_green}{RGB}{145, 172, 107}
\definecolor{acc_lightblue}{RGB}{ 65, 190, 255}
\definecolor{acc_lightgreen}{RGB}{181, 202, 130}
\definecolor{acc_orange}{RGB}{255, 128,   0}
\definecolor{acc_red}{RGB}{229,  52,  24}
\definecolor{acc_yellow}{RGB}{255, 180,   0}

\definecolor{diag_darkgreen}{RGB}{  0, 124,  48}
\definecolor{diag_darkred}{RGB}{156,  13,  22}
\definecolor{diag_gold}{RGB}{249, 186,   0}
\definecolor{diag_green}{RGB}{103, 154,  29}
\definecolor{diag_orange}{RGB}{214,  76,  13}
\definecolor{diag_pantone283}{RGB}{152, 198, 234}
\definecolor{diag_pantone300}{RGB}{  0, 101, 189}
\definecolor{diag_pantone301}{RGB}{  0,  82, 147}
\definecolor{diag_pantone540}{RGB}{  0,  51,  89}
\definecolor{diag_pantone542}{RGB}{100, 160, 200}
\definecolor{diag_petrol}{RGB}{  0, 119, 138}
\definecolor{diag_purple}{RGB}{105,   8,  90}
\definecolor{diag_red}{RGB}{196,   7,  27}
\definecolor{diag_tumgreen}{RGB}{162, 173,   0}
\definecolor{diag_tumivory}{RGB}{218, 215, 203}
\definecolor{diag_tumorange}{RGB}{227, 114,  34}
\definecolor{diag_violet}{RGB}{ 15,  27,  95}
\definecolor{diag_yellow}{RGB}{255, 220,   0}

\colorlet{diag_darkgreen_85}{diag_darkgreen!85!white}
\colorlet{diag_darkgreen_70}{diag_darkgreen!70!white}
\colorlet{diag_darkgreen_55}{diag_darkgreen!55!white}
\colorlet{diag_darkred_85}{diag_darkred!85!white}
\colorlet{diag_darkred_70}{diag_darkred!70!white}
\colorlet{diag_darkred_55}{diag_darkred!55!white}
\colorlet{diag_gold_85}{diag_gold!85!white}
\colorlet{diag_gold_70}{diag_gold!70!white}
\colorlet{diag_gold_55}{diag_gold!55!white}
\colorlet{diag_green_85}{diag_green!85!white}
\colorlet{diag_green_70}{diag_green!70!white}
\colorlet{diag_green_55}{diag_green!55!white}
\colorlet{diag_orange_85}{diag_orange!85!white}
\colorlet{diag_orange_70}{diag_orange!70!white}
\colorlet{diag_orange_55}{diag_orange!55!white}
\colorlet{diag_pantone283_85}{diag_pantone283!85!white}
\colorlet{diag_pantone283_70}{diag_pantone283!70!white}
\colorlet{diag_pantone283_55}{diag_pantone283!55!white}
\colorlet{diag_pantone300_85}{diag_pantone300!85!white}
\colorlet{diag_pantone300_70}{diag_pantone300!70!white}
\colorlet{diag_pantone300_55}{diag_pantone300!55!white}
\colorlet{diag_pantone301_85}{diag_pantone301!85!white}
\colorlet{diag_pantone301_70}{diag_pantone301!70!white}
\colorlet{diag_pantone301_55}{diag_pantone301!55!white}
\colorlet{diag_pantone540_85}{diag_pantone540!85!white}
\colorlet{diag_pantone540_70}{diag_pantone540!70!white}
\colorlet{diag_pantone540_55}{diag_pantone540!55!white}
\colorlet{diag_pantone542_85}{diag_pantone542!85!white}
\colorlet{diag_pantone542_70}{diag_pantone542!70!white}
\colorlet{diag_pantone542_55}{diag_pantone542!55!white}
\colorlet{diag_petrol_85}{diag_petrol!85!white}
\colorlet{diag_petrol_70}{diag_petrol!70!white}
\colorlet{diag_petrol_55}{diag_petrol!55!white}
\colorlet{diag_purple_85}{diag_purple!85!white}
\colorlet{diag_purple_70}{diag_purple!70!white}
\colorlet{diag_purple_55}{diag_purple!55!white}
\colorlet{diag_red_85}{diag_red!85!white}
\colorlet{diag_red_70}{diag_red!70!white}
\colorlet{diag_red_55}{diag_red!55!white}
\colorlet{diag_tumgreen_85}{diag_tumgreen!85!white}
\colorlet{diag_tumgreen_70}{diag_tumgreen!70!white}
\colorlet{diag_tumgreen_55}{diag_tumgreen!55!white}
\colorlet{diag_tumivory_85}{diag_tumivory!85!white}
\colorlet{diag_tumivory_70}{diag_tumivory!70!white}
\colorlet{diag_tumivory_55}{diag_tumivory!55!white}
\colorlet{diag_tumorange_85}{diag_tumorange!85!white}
\colorlet{diag_tumorange_70}{diag_tumorange!70!white}
\colorlet{diag_tumorange_55}{diag_tumorange!55!white}
\colorlet{diag_violet_85}{diag_violet!85!white}
\colorlet{diag_violet_70}{diag_violet!70!white}
\colorlet{diag_violet_55}{diag_violet!55!white}
\colorlet{diag_yellow_85}{diag_yellow!85!white}
\colorlet{diag_yellow_70}{diag_yellow!70!white}
\colorlet{diag_yellow_55}{diag_yellow!55!white}